\documentclass[a4paper,11pt]{article}
\usepackage[utf8]{inputenc}
\usepackage[T1]{fontenc}
\usepackage{graphicx}
\usepackage{url} 
\usepackage{amsmath, amsthm, amssymb}
\usepackage[hidelinks]{hyperref} 
\usepackage[shortlabels]{enumitem} 
\usepackage{tikz}
\usepackage{array}
\usepackage[font=small, labelfont={sf,bf}, margin=1cm]{caption}
\newcommand{\condprob}[2]{\raise2pt
                   \hbox{%
                   \mathsurround=0pt$#1$}
                    \ \raise-1pt\hbox{\scalebox{1.2}[1.5]{/}}\,%
                   \raise-2pt
                   \hbox{%
                   \mathsurround=0pt$#2$}
                    }
\newcommand{\votes}[4]{
{\arraycolsep=1.4pt\def\arraystretch{0.9}
\footnotesize
\ensuremath
\begin{matrix}  #1 &  #2 \\ #3 & #4\end{matrix}}
}                   
\renewcommand{\binom}[2]{\frac{n!}{x!y!z!t!}}

\usepackage[margin=1in]{geometry}
    \setlist[itemize]{topsep=0pt, partopsep=0pt, itemsep=\parskip, parsep=0pt}
    \setlist[enumerate]{topsep=0pt, partopsep=0pt, itemsep=\parskip, parsep=0pt}
    \setlist[description]{topsep=0pt, partopsep=0pt, itemsep=\parskip, parsep=0pt}
\usepackage{fancyhdr}
\theoremstyle{plain}
\newtheorem{theorem}{Theorem}[section]
\newtheorem{proposition}[theorem]{Proposition} 
\newtheorem{lemma}[theorem]{Lemma}
\newtheorem{corollary}[theorem]{Corollary}
\theoremstyle{definition}    
\newtheorem{definition}[theorem]{Definition}
\newtheorem{example}[theorem]{Example}

\begin{document}
\title{{Optimal decision rules for the discursive dilemma}
}
\author{Aureli Alabert \\ 
           Department of Mathematics \\
           Universitat Aut\`onoma de Barcelona \\  
           08193 Bellaterra, Catalonia \\  
           \url{Aureli.Alabert@uab.cat}  
           \and
           Mercè Farré \\ 
           Department of Mathematics \\
           Universitat Aut\`onoma de Barcelona \\  
           08193 Bellaterra, Catalonia \\  
           \url{farre@mat.uab.cat}
           \and
           Rubén Montes \\
           Department of Mathematics \\
           Universitat Aut\`onoma de Barcelona \\  
           08193 Bellaterra, Catalonia \\            
}
\maketitle
\begin{abstract}
We study the classical \emph{discursive dilemma} from the point of view of finding the best 
decision rule according to a quantitative criterion, under very mild restrictions on the
set of admissible rules. The members of the deciding committee are assumed to have
a certain probability to assess correctly the truth or falsity of the premisses, 
and the best rule is the one that minimises a combination of the probabilities
of false positives and false negatives on the conclusion. 
\end{abstract}

{\small
\textbf{Keywords}: Discursive dilemma, doctrinal paradox, judgment aggregation, truth tracking. 
}

\section{Introduction}\label{sec:Intro}
\subsection{Statement of the problem}
Nowadays, the so-called \emph{doctrinal paradox} is a classical problem in judgment aggregation,
in which different reasonable majority-type voting rules may lead to different conclusions:
A group of people must assess the simultaneous truth or not of a set of premisses, and voting
first on each premiss or voting directly on the conclusion not necessarily yield the same
result. A slightly different formulation of the paradox have received the name of 
\emph{discursive dilemma}.

In practice, the situation may appear when a court is deciding if a defendant is guilty (a set
of evidences are all verified), or not-guilty (at least one of the evidences is false), and this is the
origin of the name (Kornhauser \cite{Kornhauser1992169}). 
But obviously it is present beyond legal cases. A prize or a job position can 
be awarded if and only if several debatable conditions concur; or several subjective medical indicators
determine the presence of an illness or the need of a treatment; and so on.
As soon as three people join to make a decision on a compound question, the paradox is potentially
present.

To precise, suppose there are 
two clauses $P$ and $Q$ and each member of a committee has to decide between $P$ and its
negation $\neg P$ 
and between $Q$ and its negation $\neg Q$; and that the final goal is to assess if $C:=P\wedge Q$ is true 
or its negation $\neg C$ 
is true. In the court example, the jury has to decide if the defendant is guilty 
or not, and it is agreed beforehand that the guilty verdict 
is logically equivalent to the truth of both premisses $P$ and $Q$.

Suppose the voters first decide by simple majority between $P$ and $\neg P$, and 
separately between $Q$ and $\neg Q$. If both $P$ and $Q$ get the majority, then the 
conclusion is $C$, and otherwise it is $\neg C$. This decision rule is called 
\emph{premiss-based}. Suppose on the other hand that  
each voter decides directly on $C$ or $\neg C$, and then the collective decision
is taken by simple majority on these alternatives. This rule is called \emph{conclusion-based}.

There are cases where the premiss-based rule leads to $C$, while the conclusion-based yields $\neg C$.
For instance, for a 3-member committee, this happens when one of the members thinks 
both $P$ and $Q$ are true, the second one thinks $P$ is true and $Q$ is false, and the 
last one thinks $P$ is false and $Q$ is true. 
Sometimes it is said that the conclusion-based rule is more ``conservative'' than the premiss-based
rule (or that the latter is more ``liberal'' than the former), because
the positive conclusion $C$ is frequently the ``risky one''.

These two rules are quite natural, and both can be justified on intuitive or philosophical
grounds; see for example, Mongin \cite[section 2]{Mongin2012-MONTDP}. In particular,
the conclusion-based rule respects the deliberation of the individual judges; in the premiss-based
rule the decision can be fully justified in legal terms.
Others rules can be proposed. 
In \cite{AlabertFarre2022}, we introduced a new rule, 
which stands in some sense midway between the premiss-based and the conclusion-based
rules.

\bigskip
In this paper we study general decision rules for the situation given.
We will consider
the set of all possible decision rules, subject only to a very mild rationality requirement.
They will be called \emph{admissible rules}. 
We want to study this set as a whole, 
and find the best rule according to some objective criterion, disregarding whether that rule can or cannot 
be explained on ``logical'' or ``intuitive'' grounds, it is a consequence of some political or sociological 
idea, or it satisfies some other desirable property.

The mentioned rationality requirement states only 
that if a member of the committee changes their\footnote{The singular \emph{they/their} will be used 
to avoid gender bias.} opinion on a clause in some direction, 
the conclusion can only eventually change in the same direction.

We take the epistemic point of view
 that \emph{there is an actual truth}
that we want to guess with the highest possible confidence. This is different from 
the aggregation 
of preferences as in elections, or in taking decisions on the course of actions, where 
there is not an absolute truth. 

Our objective criterion is related to the minimisation of the combined chances to incur 
in false positives (deciding
$C$ when the reality is $\neg C$) and in false negatives (deciding $\neg C$ when the 
reality is $C$). This is explained in detail in Section \ref{sec:ModelCrit} and of course involves
a mathematical (probabilistic) setting where all elements have to be precisely defined.
Our point of view is thus ``conclusion-centric'', in the sense that we do not care about the
correct guessing of the premisses.

We emphasise that the adoption of a particular criterion is a modelling choice, 
and it is what confers the rationale to the best rule under it. 
The criterion proposed here can be replaced by another one, if deemed better for the
situation at hand, 
and the philosophy of \emph{finding the best under the chosen criterion} can be
applied as well. We consider here, in fact, a family of criteria, parametrised by the relative weight
put on false positives and false negatives. 

With this optimisation approach, we do not need to talk about \emph{majorities}. The votes of the $n$ members 
 of the committee will be split into four slots: $P\wedge Q$, $P\wedge \neg Q$,
 $\neg P\wedge Q$, and $\neg P\wedge\neg Q$, and the aggregated number of votes for each 
 possibility will be non-negative integers $x,y,z,t$, respectively, with 
 $x+y+z+t=n$, being $n$ the number of voters.
 In the case of three premisses, there will be 8 slots, and in general  
 $p$ premisses would give $2^p$ different possible votes of each member.
 A decision rule states, for each possible values of $x,y,z,t$ which decision, $C$ or $\neg C$,
 is taken.
 
The number of rules grows exponentially with $n$. 
The admissible rules are much less, and they can be implicitly enumerated 
so that all computations needed to find the optimal rule or a ranking of
rules are relatively efficient.

\bigskip
If each committee member could infallibly guess the truth or falsity
of each premiss, then the correct truth or falsity of the conclusion will be reached
without difficulty. In fact, a single-member committee would suffice. The whole point
of having multi-member committees is to alleviate the possibility that the final conclusion
be wrong. It is therefore quite natural to use a probabilistic model that starts 
with the (estimated) probability that the committee members make the correct guessing
on each premiss. We call this probability their \emph{competence}, and we assume that it is
greater that $\frac{1}{2}$, and that is the same for all members of the committee and
for all premisses, although this is easily relaxed, as we will see in the final section.

The collective decision guesses correctly or incorrectly with some probability
that depends on the voters' competence and on the real truth value of the premisses.
Only the 
conclusion matters, and only the premisses are voted. One may think, as pointed out
by Mongin \cite{Mongin2012-MONTDP}, that an external judge has to decide on the conclusions
after the committee has sent them their individual opinions.

\subsection{Related literature}

\textbf{Doctrinal paradox.}
The term \emph{doctrinal paradox} appears first in the works of Kornhauser \cite{Kornhauser1992b}, 
and Kornhauser and Sager \cite{KornhauserSager1993}.
They were interested in legal court cases, so that they
spoke of \emph{issue-by-issue} and \emph{case-by-case} majority voting.

Pettit \cite{Pettit2001} and List and Pettit \cite{list_pettit_2002} formulate the problem
in terms of propositional logic, and called it the \emph{discursive dilemma}. The simple example
of the three-member committee cited above, can be summarised in Table \ref{tab:dilemma}:
\begin{table}[h] \centering
    {\extrarowheight 2pt
  \begin{tabular}{r|c|c|c}
     Voter & Proposition $P$ & Proposition $Q$ & Proposition $C$ \\
    \hline 
    $1$ & $\phantom{\neg}P$ & $\phantom{\neg}Q$ & $\phantom{\neg}C$ \\
    $2$ & $\phantom{\neg}P$ & $\neg Q$ & $\neg C$ \\
    $3$ & $\neg P$ & $\phantom{\neg}Q$ & $\neg C$ \\
    \hline 
    majority & $\phantom{\neg}P$ & $\phantom{\neg}Q$ & $\neg C$ \\
  \end{tabular}
}
\caption{The discursive dilemma: The collective majority voting in the three premisses
is inconsistent with the doctrine $C\Leftrightarrow P\wedge Q$.}\label{tab:dilemma}
\end{table} 
  In the Kornhauser--Sager formulation, the committee votes either on the first
  two propositions (premiss-based/issue-by-issue), or on the third (conclusion-based/case-by-case),
  and the two results are different. In the List--Pettit formulation, the 
  committee votes on the three propositions, and this leads to a logical inconsistency. 
  The inconsistency comes from the constraint $C\Leftrightarrow P\wedge Q$
  (the ``doctrine'' to which there is a previous agreement). We see that the individual members of
  the committee adhere to the doctrine; however the committee as a whole does not.
  
  The advantage of the formulation in terms of propositional logic is that it can be
  generalised to any set of propositions, to the point that the distinction between premisses and 
  conclusions may be unnecessary. In general, an \emph{agenda} is a logically consistent set of
  propositions, closed under negation, on which judgments have to be made, and that can be 
  entangled by  logical constraints. In our case the agenda is 
  $\big\{P,\neg P, Q,\neg Q, C, \neg C\big\}$, with the constraint $C\Leftrightarrow P\wedge Q$. 
  In this setting, the 
  doctrinal paradox (or more properly, the discursive dilemma) reads: 
  If all pairs of formulae in the agenda are decided by majority, the resulting set
  of propositions can be inconsistent. 

\textbf{Judgment aggregation.}
  The body of knowledge that has been developed from List--Pettit formulation is known as Judgment Aggregation
  Theory (or Logical Aggregation Theory, as proposed by Mongin \cite{Mongin2012-MONTDP}). 
  In a quite natural way, the backbone of the theory is formed by (im)possibility results
  on the existence of aggregation rules satisfying certain desirable axioms. List and Pettit
  \cite{list_pettit_2002}, \cite{ListPettit2004} already proved results of this kind, 
  extended very soon by  Pauly and van Hees \cite{PaulyVanHess2006}, Dietrich \cite{DIETRICH2006286},
  and Nehring and Puppe \cite{NehringPuppe2008}.
  
  The aggregation problem is described in full generality 
  for example in the preliminaries of Nehring and Pivato \cite{NEHRING20111488}
  and Lang et al.~\cite{Lang2017}, and in the complete surveys by 
   Mongin \cite{Mongin2012-MONTDP}, List and Puppe \cite{List2009-LISJAA}, List and Polak \cite{ListPolak2010441}, 
   and List \cite{List2012}.
  A \emph{judgment} is defined as a mapping
  from the agenda to the doubleton $\{\text{True}, \text{False}\}$; 
  a \emph{feasible judgment} respects moreover the  
  underlying logical constraints of the propositions\footnote{Lang et al.~\cite{Lang2017} call it 
    \emph{consistent judgment}, in the sense that it is logically consistent
    (not a contradiction) when the logical constraints are added.}. The \emph{judgment 
  aggregation problem} is then defined as the construction of a feasible reasonable 
  collective judgment from the voters' individual judgments. Formally, an 
  \emph{aggregation rule} $F$ is a mapping that assigns to every \emph{profile}
  $(J_1,\dots,J_n)$ of individual judgments $J_i$ of the $n$ voters, 
  a collective judgment $J=F(J_1,\dots,J_n)$. A \emph{feasible aggregation rule}
  must assign a feasible judgment to any input of feasible profiles. 
  Feasible rules trivially exist; for instance $J\equiv J_i$ for some $i$
  is such a rule (called a \emph{dictatorship}, for obvious reasons). 
  Banning dictatorships and imposing other mild desirable conditions 
  leads very quickly to non-existence of feasible aggregation rules.
  The range of possible \emph{voting paradoxes} is the set
  of non-feasible mappings. The classical doctrinal paradox case, despite its simplicity, already features one such non-feasible
  mapping, namely $P\mapsto \text{True}$, $Q\mapsto \text{True}$, $C\mapsto \text{False}$.
  
  In a \emph{truth-functional} agenda, the propositions
  are split into a set of premisses, and a set of conclusions. Assigning a Boolean value to all
  premisses, and applying the logical constraints, the value of all conclusions is determined.
  This is clearly the case in the doctrinal paradox setting, where moreover the premisses consist 
  of mutually independent (not linked by constraints) proposition-negation pairs. 
  The truth-functional case in general has been studied mainly in 
  Nehring and Puppe \cite{NehringPuppe2008}, for independent as well as interdependent premisses,
  and in Dokow and Holzman \cite{DokowHolzman2009} (see also Miller and Osherson \cite{MillerOsherson2009}). 

\textbf{Distance-based methods and truth-tracking.}
    From 2006 (Pigozzi \cite{Pigozzi2006}, Dietrich and List \cite{DietrichList2007}) 
    another point of view emerged, in 
    which specific judgment rules are proposed, and their properties studied. See Lang et al.~\cite{Lang2017}
    for a partial survey, and the references therein. 
    Most of these rules can be defined as some sort of optimisation with respect to a criterion, 
    i.e.~the rule is defined as the one(s) that maximises or minimises a certain quantity,
    usually a distance or pseudo-distance to the individual profiles, while providing a 
    consistent consensus judgment set. 
 
            
There are two different
    situations to which they can be applied. Either the collective judgment set is a decision on 
    the course of actions (as in the adoption of public policies), or there is an underlying 
    objective truth of each proposition under scrutiny 
    that one would like to guess (as in court cases). 
    The latter is called \emph{truth-tracking} (or \emph{epistemic}) judgment aggregation,
    and it is where the present work belongs. 
 
  Thus the goal is to get the right values
  of this pre-existent ``state of nature'', or at least the right values of the set of conclusions
  in the truth functional case.
  In this context, the concept of \emph{competence} 
  of the voters arise naturally: 
  how likely is that a voter guess the correct answer to an issue? And it is also
  natural to model this likelihood as a probability. Actually, this approach dates back to 
  Condorcet and his celebrated Jury Theorem.

  The competence as a parameter has been studied, for example, in 
  Bovens and Rabinowicz \cite{Bovens2006} and in Grofman et al.~\cite{Grofman1983} for the
  one-issue case. The latter extends the Condorcet theorem in several directions, 
  particularly for the case of unequal competences among voters.    
        
  List \cite{List2005} computed the probability of appearance of the doctrinal paradox,
  and the probability of correct 
  truth-tracking as a function of the different states of nature, allowing for different competences 
  in judging both premisses but the same across individuals.
  Fallis \cite{Fallis2005} also observed that the premiss-based rule is better or not than the conclusion-based rule
  depending on the competence and on the ``scenario'' (state of nature). 
  
  The fact that the probability to guess correctly the truth depend on the unknown
  true state of nature leads to a modelling choice: Either we specify an a priori 
  probability distribution on the possible states of nature (for the doctrinal paradox,
  the four states $P\wedge Q$, $P\wedge\neg Q$, $\neg P\wedge Q$ and $\neg P\wedge\neg Q$),
  or we have to resort to conservative estimations, as in classical (non-Bayesian) 
  statistics. The main approach in this paper is the second, but most of the related
  literature assumes the first. Notably:  
  
  The cited paper by Bovens and Rabinowicz \cite{Bovens2006} compares the premiss-based
  and the conclusion-based rules under the assumption of same competence for 
  both premisses and their negations, and independence of voters, as in the present paper.
  They impose a Bernoulli prior on each premiss, the same for all.
   
  Hartmann et al.~\cite{Hartmann-al2010} 
  aims at generalising \cite{List2005} and \cite{Bovens2006}
  with a conjunctive truth-functional agenda allowing more
  than two premisses. The authors propose a continuum of distance-based rules, parametrised by the weight 
  of the conclusion relative to the premisses, 
  and containing the premiss- and conclusion-based procedures as extreme cases. 
  The hypotheses are essentially the same as in \cite{Bovens2006}. 
  Miller and Osherson \cite{MillerOsherson2009} also propose a variety of distance metrics, and distinguish
  between ``underlying metric'' and ``solution method''. Each solution method chooses a loss function to
  minimise (based on the metric) and a set of eligible rules.
    
  The point of view of Pivato \cite{Pivato2013} is that the votes are observations of the `truth plus noise'.
  This allows to think of the profile of individual judgments as a statistical sample (at least under
  the hypothesis of the same noise distribution for all voters), and study the decision rules as
  statistical estimators.

\textbf{Truth-functional agendas and truth-tracking.}
  The abovementioned papers List \cite{List2005}, Fallis \cite{Fallis2005}, Bovens and Rabinowicz \cite{Bovens2006}, 
  Hartman et al.~\cite{Hartmann-al2010} and, partially, Miller and Osheron \cite{MillerOsherson2009}, 
  deal with truth-functional agendas.
  Nehring and Puppe \cite{NehringPuppe2008} focusses  on truth-functional agendas with logically interdependent premisses.
   
  Combined with a truth-functional agenda, the truth-tracking setting can still be concerned either with 
  guessing the truth of all propositions (`getting the right answer for the right reasons')
  or only on the conclusion (`getting the right answer for whatever reasons').
  Bovens and Rabinowicz \cite{Bovens2006} discuss the merits of premiss- and conclusion-based procedures
  for both goals. 
  Bozbay et al.~\cite{BozbayDietrichPeters2014} and Bozbay \cite{Bozbay2019} also study both aims, for 
  independent and for interrelated issues, respectively. The cited work by Hartmann et al.~\cite{Hartmann-al2010} 
  is conclusion-centric (``whatever reasons''), while Pigozzi et al.~\cite{Pigozzi2009}, being conclusion-centric, 
  applies later a procedure based on Bayesian networks to get the premisses that ``interpret'' the
  previously decided conclusions.

  Distance-based methods are nothing else that the minimisation of a \emph{loss function} that measures
  the dissatisfaction with every possible consistent outcome. Equivalently, one may maximise \emph{utility
  functions}. Both are capable to account for the consequences of the decisions, and thus allow to 
  set up more complete models, in line with Statistical Decision Theory (Berger \cite{Berger1985}).
  Different loss functions or utilities give rise to possibly different optimal rules, and 
  it is a modelling task to choose the right loss function for the problem at hand. 

  In this sense, Fallis \cite{Fallis2005} writes about the `epistemic value', highlighting that
  guessing correctly a proposition may have a different value that guessing correctly its contrary;
  Bozbay \cite{Bozbay2019} uses a simple 0-1 utility function to indicate incorrect-correct
  guessing (of all propositions or of the conclusions alone); 
  Hartmann et al.~\cite{Hartmann-al2010} tries giving different utilities to false positives and false
  negatives on the conclusion to assess the performance of their continuum of metrics; 
  finally, Bovens and Rabinowicz \cite{Bovens2006}, in the discussion section, suggest  
  introducing different utilities to each correct guessing to compare the premiss-based and the 
  conclusion-based voting rules in each practical case.

\bigskip
Our proposal in this paper is an optimisation criterion in the truth-tracking, conclusion-centric
    case of the discursive dilemma. The best decision rule minimises a combination of 
    false positives and false negatives,
    and any two rules can be easily compared according to this criterion.
    No prior on the states of nature needs to be established, although it can be easily accommodated. 
    In the theoretical results we assume the same competence level of all committee members and all premisses,
    and independence among premisses. In practice, the specific computation of the score of each rule can
    be done under much less assumptions. 
    In any case, the loss function fully determines the optimal rule.
  
  We do not consider strategic voting; we assume
  that everyone votes honestly each of the premisses. Strategic voting is 
  conceivable even in the simple doctrinal paradox case: someone who honestly would vote for 
  $P$ and $\neg Q$ could change to $\neg P$ and $\neg Q$ just to push more
  for the $\neg C$ conclusion). Strategic voting is considered in 
   Bozbay et al.~\cite{BozbayDietrichPeters2014}, de Clippel and Eliatz \cite{deClippel201534},
   Terzopolou and Endriss \cite{TerzopoulouEndrissSAGT2019} and Bozbay \cite{Bozbay2019}.

\subsection{Organisation of the paper}
The remainder of the paper is organised as follows: In Section \ref{sec:DecRules} we study the structure of the sets of 
voting tables and of decision rules, which are both partially ordered sets with an order induced by the 
admissibility condition. 

Section \ref{sec:ModelCrit} introduces in the first part our probabilistic model, based on the probability 
that each committee member guesses correctly the truth value of each premiss, 
and the concepts of false positive and false negative when the true state of nature 
is unknown. 
In the second part, we introduce the family of optimisation
criteria, parametrised by a relative weight assigned to the two errors. 

Section \ref{sec:Main} contain our main theoretical results. It turns out that it is relatively easy to determine whether
a given voting table leads to the conclusion $C$ or $\neg C$ in the optimal rule. Moreover, this depends
only on two numbers: the difference between votes to $P\wedge Q$ and to $\neg P\wedge \neg Q$, and 
the difference between votes to $P\wedge \neg Q$ and to $\neg P\wedge Q$. This simplifies the structure
of the set of voting tables, and shortens the evaluation of the rules. In this section we also
characterise completely the set of values of competence and weight for which the premiss-based
rule is optimal. 

Section \ref{sec:ComputSoft} explains details on the actual computations, and describes the accompanying software,
downloadable from \url{https://discursive-dilemma.sourceforge.io}.
Finally, in Section \ref{sec:ConcDisc} we discuss the 
results and propose possible extensions. Some marginal computations and checks have been left to an appendix.
\section{The set of decision rules}\label{sec:DecRules}
A possible voting result of a committee with $n$ members assessing on issues $P$ and $Q$ will be a 
\emph{table}, denoted  
$\Big[\,\votes{x}{y}{z}{t}\,\Big]$,
or $(x,y,z,t)$ to save space, with non-negative integer entries, adding up to $n$,
representing the quantity of votes received by the 
options $P\wedge Q$, $P\wedge \neg Q$, $\neg P\wedge Q$ and $\neg P\wedge\neg Q$, respectively. 
The set of all such tables will be denoted by $\mathbb T$. 
 
 A \emph{decision rule} can be thought of as a mapping $\mathbb T\longrightarrow \{0,1\}$. 
 Tables mapped to 
 $1$ are those that entail the decision $C=P\wedge Q$; those mapped to $0$ represent the 
 opposite, $\neg C$. Sometimes we will call them \emph{positive tables} and \emph{null tables}, 
 and denote by $\mathbb T^+$ and $\mathbb T^0$ the respective sets. 
 The decision rule can also be seen as the subset of positive tables,  
 and we will make use of both interpretations. 
 
 There are $N=\binom{n+3}{3}$ ways to fill a voting table, and $2^N$ 
 decision rules, as many as subsets of the set of tables. This is a huge number, already for $n=3$, 
 but the set of ``reasonable rules'' will be much more modest. 

 Two tables $(x,y,z,t)$ and $(x,z,y,t)$ are \emph{transposed} of each other. Since $P$ and
 $Q$ will play symmetric roles, it makes sense to admit only rules that assign the same decision
 to both tables. 
 
\begin{figure}[t]
\centering
\begin{tikzpicture}
    \node [anchor=south west,inner sep=0] at (0,-0.55) {\includegraphics[scale=0.45]{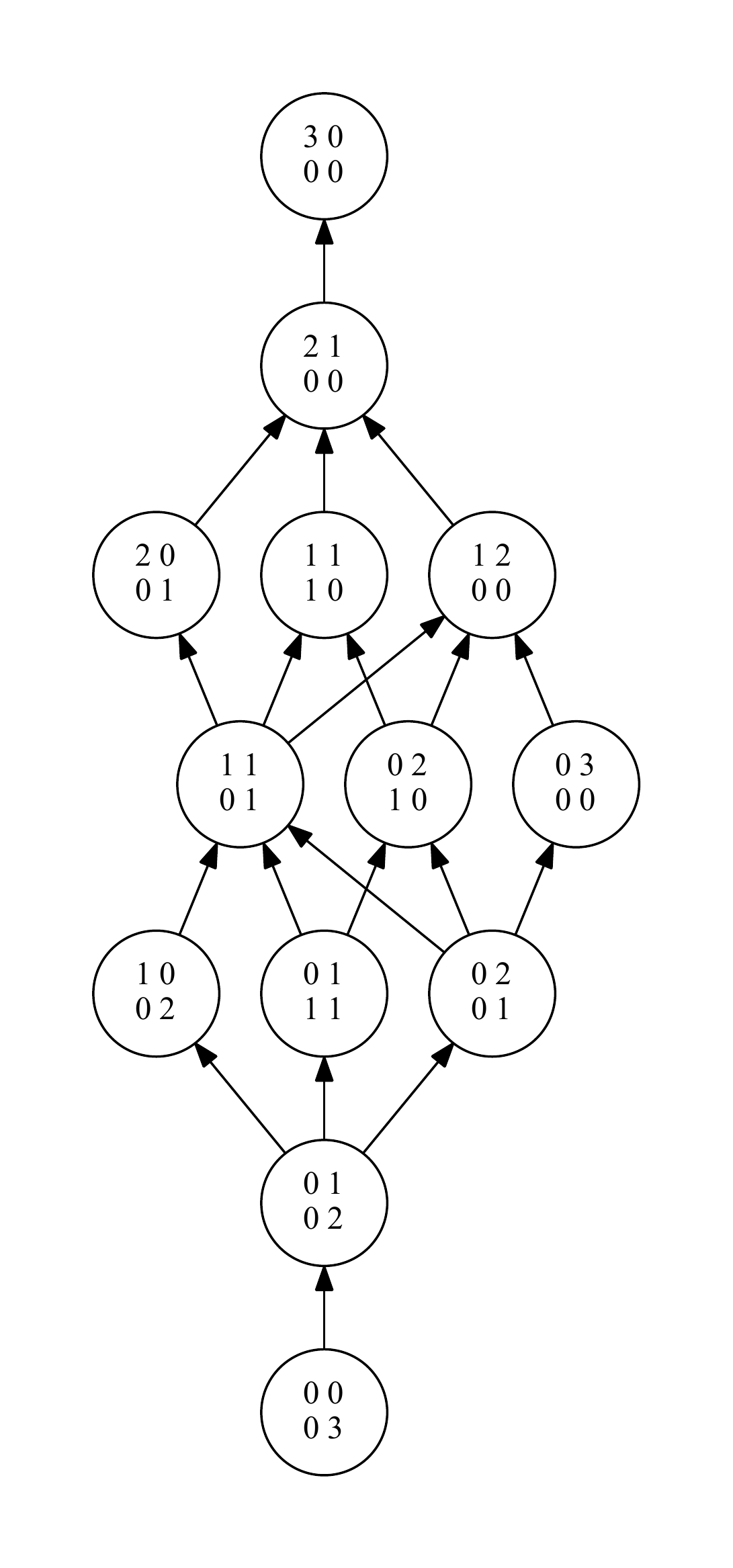}};
        \node [left] at (-0.3, 0.5) {\small $\rho=-3$};
        \node [left] at (-0.3, 1.93) {\small $\rho=-2$};
        \node [left] at (-0.3, 3.36) {\small $\rho=-1$};
        \node [left] at (-0.3, 4.80) {\small $\rho=0$};
        \node [left] at (-0.3, 6.22) {\small $\rho=1$};
        \node [left] at (-0.3, 7.65) {\small $\rho=2$};
        \node [left] at (-0.3, 9.08) {\small $\rho=3$};
\end{tikzpicture}
\caption{Hasse diagram of the poset $(\mathbb T,\le)$ for committee size $n=3$, 
and transposed tables identified. The arrows correspond to the
transitive reduction. All other relations are deduced by transitivity.
The rank function $\rho(x,y,z,t)=x-t$ is also represented.}
\label{fig:n3graph}
\end{figure}
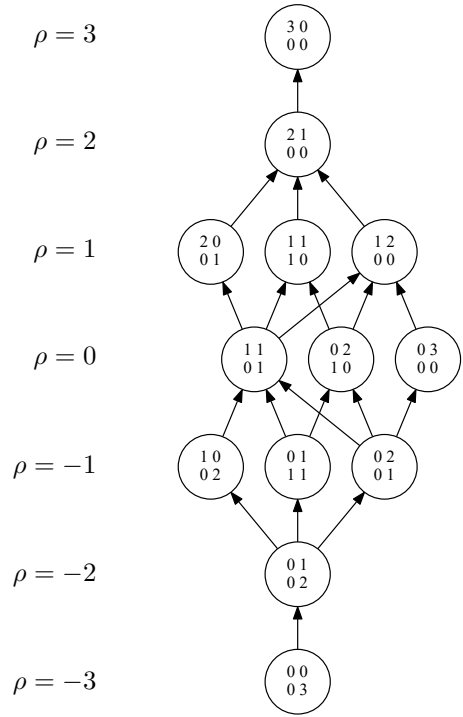

 Besides this symmetry, we impose another condition for the admissibility of a 
 decision rule. Suppose that, given
 a positive voting table $(x,y,z,t)$, one of the voters of $\neg P$ changes their choice to $P$,
 or one of the voters of $\neg Q$ changes to $Q$. The new table would support better
 the conclusion $C$ than the older; hence it makes sense to impose that the new table be
 also a positive table. To formulate the condition in a
 mathematically practical way, let us introduce the partial order on $\mathbb T$
 generated by the four relations
\begin{equation}\label{eq:transreduc}
\begin{split}
(x,y,z,t)&\le(x,y,z+1,t-1) \\
(x,y,z,t)&\le(x,y+1,z,t-1) \\
(x,y,z,t)&\le(x+1,y-1,z,t) \\
(x,y,z,t)&\le(x+1,y,z-1,t) \ ,
\end{split} 
\end{equation}
that means, the smallest partial order $\le$ that satisfies relations (\ref{eq:transreduc}) above, for all $x,y,z,t$ for
which they make sense. A partially ordered set is also called a \emph{poset}, for short.
The relations (\ref{eq:transreduc}) are called the \emph{transitive reduction}
of the poset $(\mathbb T,\le)$. Posets can be represented by \emph{Hasse diagrams}, 
which are directed graphs with the transitive reduction represented by arrows pointing in the increasing direction.
The case of committee size $n=3$ is depicted in 
Figure \ref{fig:n3graph}, where transposed tables have been identified; they occupy the same
spot and are not comparable. 
When two tables $T,S\in\mathbb{T}$, satisfy $T\le S$ and $T\neq S$, we shall obviously write $T<S$, 
or $S>T$.

We thus arrive to the following reasonable definition of admissibility:
\begin{definition}\label{def:admissible}
A decision rule $r\colon \mathbb T\rightarrow \{0,1\}$ is \emph{admissible} if:
\begin{enumerate}
\item It takes the same value on transposed tables:
\begin{equation*}
r(x,y,z,t)=r(x,z,y,t) \ .
\end{equation*}
\item It is order-preserving on the partially ordered set $(\mathbb T,\le)$:
\begin{equation*}
(x,y,z,t)\le(x',y',z',t') \Rightarrow r(x,y,z,t)\le r(x',y',z',t') 
\ .
\end{equation*}
\end{enumerate} 
\end{definition}    

\begin{example}\label{ex:classicalRules}
The classical premiss-based rule $r_{pb}$ is defined by 
\begin{equation*}
\text{$r_{pb}(x,y,z,t)=1$ if and only if $x+y>z+t$ and $x+z>y+t$\ ,}
\end{equation*}
whereas the conclusion-based rule $r_{cb}$ is given by
\begin{equation*}
 \text{$r_{cb}(x,y,z,t)=1$ if and only if $x>y+z+t$\ .}
\end{equation*}
It is readily checked that both rules are admissible in the sense of Definition \ref{def:admissible}. 
In Alabert and Farré \cite{AlabertFarre2022}, another admissible rule was introduced, 
called \emph{path-based}, 
and defined by 
\begin{equation*}
 \text{\rm \emph{$r_{hb}(x,y,z,t)=1$ if and only if $x>z+t$ and $x>y+t$\ .}}
\end{equation*}
As an example of a non-admissible rule, consider declaring $C$ true if and only if the votes for 
$P\wedge Q$ are more than any other combination of premisses (i.e. $x>y$ and $x>z$ and $x>t$).
\qed
\end{example} 
We will see later that the second admissibility condition is not restrictive with respect
to our optimisation criterion: given a non-order-preserving rule, there exists 
an order-preserving one that performs better.

The poset $(\mathbb T,\le)$ is \emph{ranked} (also called \emph{graded}), i.e.~there exists 
a \emph{rank function} $\rho$ compatible with the order relation:
It satisfies  
$T<S\Rightarrow\rho(T)<\rho(S)$, 
and if 
$S$ is an \emph{immediate successor} of T (there are no elements in between),
then $\rho(S)=\rho(T)+1$. 
In the Hasse diagram,
each rank can be pictured as a ``level'' in the graph (see again Figure \ref{fig:n3graph}).

To prove that $(\mathbb T, \le)$ is a ranked poset, we use the known result that a finite poset
admits a rank function if and only if all maximal chains have
the same length.\footnote{Called the \emph{Jordan--Dedekind property}.} Recall that a \emph{chain} is 
a totally ordered subset of the poset. 
A \emph{maximal chain} is a chain with maximal cardinality.

A poset is \emph{connected} if for every two elements T and S there is a 
finite sequence $T=U_1,\dots,U_n=S$ of elements  
such that $U_i$ and $U_{i+1}$ are comparable, i.e. 
either $U_i\le U_{i+1}$ or $U_{i+1}\le U_{i}$. 
The poset $(\mathbb{T},\le)$ is connected, because we can transform a table into any other
one by moving votes one at a time through the transitive reduction.

\begin{proposition}
$(\mathbb T, \le)$ is a ranked poset, with $\rho(x,y,z,t)=x-t$ as a rank function.
\end{proposition}

\begin{proof}
There is a unique minimal element, namely $m=(0,0,0,n)$, and a unique maximal 
element $M=(n,0,0,0)$. Since $(\mathbb{T},\le)$ is connected, all maximal chains
start and finish in these
elements. To go from $m$ to $M$, each vote must make two steps, one of them up or to the left
of the table (that means, from $t$ to $y$ or $z$), and the other one to the left or
up respectively (from $y$ or $z$ to $x$). These individual movements are the 
transitive reduction of the partial order $\le$,
and therefore there are no other tables in between. Since there are $n$ votes,
we need $2n$ steps to move all votes from the minimal to the maximal element, and in consequence
any maximal chain has exactly $2n+1$ elements. 

Notice that a movement towards an immediate successor imply subtracting one unit to $t$
or adding one unit to $x$, but not both. Therefore $\rho(T)=x-t$ is a rank function for
$(\mathbb T,\le)$. 
\end{proof}

The rank function of a ranked poset is not unique, but it is completely determined by 
setting the rank of any element of the poset.

\bigskip
The set $\cal A$ of 
admissible decision rules possesses also a natural partial order: 
$r\le s$ if $r(T)\le s(T)$ for all tables $T\in\mathbb T$. This is the usual partial order in
a set of real functions on any domain. Since the range of decision rules mappings is $\{0,1\}$,
the relation $r\le s$ means that the set of 
positive tables relative
to $r$ is included in the set of positive tables relative to $s$. It can be said that $r$ is less
liberal (or more conservative) than $s$ in the sense explained in the Introduction.
In terms of the risk of opting for $C$ when it is wrong, $\le$ is the relation 
``to be less risky or equal to''.

\begin{example}
Refer to the rules of Example \ref{ex:classicalRules}. The premiss-based rule is more liberal than
the path-based rule, and this one is in turn more liberal than the conclusion-based rule.
In other words, 
we are saying that, in $({\cal A},\le)$,
\begin{equation*}
r_{cb}\le r_{hb}\le r_{pb} 
\ .
\end{equation*}
This can be seen using the characterisations given in Example \ref{ex:classicalRules}
(see \cite[Proposition A.2]{AlabertFarre2022} for the proof). 
To precise a little more, one can check that rules $r_{hb}$ and $r_{cb}$ coincide for committee sizes
$n=3,5$, and they are different for $n\ge 7$. Rule $r_{pb}$ is always strictly greater than $r_{hb}$. 
\qed
\end{example} 

\bigskip
An \emph{upper set} is a subset $U$ of a poset such that 
$x\in U, x<y \Rightarrow y\in U$. 
It is immediate to see, from the second condition of admissibility, the following equivalence.
\begin{proposition}\label{pr:goodUpper}
A decision rule $r\colon {\mathbb T}\rightarrow \{0,1\}$ is order-preserving  if and only if
$\{T\in\mathbb T:\ r(T)=1\}$ is an upper set of $(\mathbb T,\le)$.
\end{proposition}

Given a poset, the family of its upper sets, with the inclusion order relation,
is a \emph{complete lattice}: A partially ordered set in which all subsets
have a \emph{supremum} (a least upper bound) and an \emph{infimum} (a greatest lower bound). 
Applied to our case, we are saying that the union and the intersection
of upper sets are upper sets or,  in terms of rule mappings, that 
the maximum ($=$~sum) and the minimum ($=$~product) of admissible rules are admissible rules. 

An \emph{antichain} is a subset $\mathbb S$ of a poset such that any two elements of $\mathbb S$ are 
not comparable. Antichains and upper sets are related in the following way: The minimal elements
of any upper set form and antichain; conversely, any
antichain $A$ determines the upper set
\begin{equation*}
  \{x:\ (\exists y)(y\in A \wedge x\ge y)\}
  \ .
\end{equation*}
The empty antichain is also considered, and corresponds in our case to the rule $r\equiv 0$.

For finite posets, the correspondence between antichains and upper sets is bijective.
Enumerating upper sets is therefore
equivalent to enumerating antichains. Even computing the number of upper sets is
not easy in general. 
For example,
in the well-known poset of the subsets of a given set of $k$ elements, 
with the inclusion relation, the 
number of upper sets (called the \emph{Dedekind numbers}, see \cite{oeis-org}),
is not known for $k>8$.

\section{Probabilistic model and optimisation criterion}\label{sec:ModelCrit}
We want to find the best of all admissible rules, according to some quantitative
criterion, formulated in terms of a probabilistic model. 
In this section the criterion will be introduced, and the next one 
will be devoted to the characterisation of the optimal rule. 

\subsection{Probabilistic model}
Suppose $C=P\wedge Q$ is the true state of nature. If for some rule $r$ and a table of votes
$T=(x,y,z,t)$, we have $r(T)=1$, we say that this is a \emph{true positive} (TP). Otherwise,
if $r(T)=0$, it is a \emph{false negative} (FN). Similarly, if $\neg C$ is the true state, 
$r(T)=0$ will be a \emph{true negative} (TN) and $r(T)=1$ will be a \emph{false positive} (FP).
Ideally, a good decision rule should minimise somehow the occurrence of false positives and false negatives.
To assess the likelihood of these occurrences we need a probabilistic model in which to evaluate the
probability of appearance of FP and FN. To that end, we need an estimate of the probability that
the members of the committee guess correctly the true value of the premisses $P$ and $Q$.

The probability that a committee member vote the correct value of the premisses will be called 
its \emph{competence}.
We will assume that all committee members have the same competence,
a number strictly between $\frac{1}{2}$ and $1$. Notice that a competence
less than $\frac{1}{2}$ does not make sense, because in that case we can 
reverse all opinions of the committee, and we get another committee with competence
greater than $\frac{1}{2}$. If it were exactly $\frac{1}{2}$, there is a trivial
solution that will be pointed out later; if it is 1, then a one-member committee is enough
and they are always right. 
We will also assume that the committee size $n$ is odd and two additional independence conditions. 
Specifically, we assume in the sequel the following
hypotheses:
\begin{description}
  \item[(H1)] 
  \emph{Odd committee size}: The number of voters is an odd number, $n=2m+1$, with $m\geq 1$.
   \item[(H2)] \emph{Equal competence}:
   The competence $\theta$ satisfies $\frac12< \theta<1$ and it is the same for all 
   voters and for both premisses $P$ and $Q$.
   \item[(H3)] \emph{Mutual independence among voters}: The decision of each voter does
   not depend on the decisions of the other voters.
   \item[(H4)] \emph{Independence between $P$ and $Q$}: For each voter, the decision on one premiss
   does not influence the decision on the other. 
\end{description}
  Formally, hypotheses (H2)--(H4) can be rephrased by saying that for each voter in the committee 
and each premiss, there is a random variable 
that takes the value 1 if the voter believes the clause is true, and zero otherwise, and 
all these random variables are stochastically independent and identically distributed. 
Their specific distribution depends on the true state of nature. 
See Section \ref{sec:ConcDisc} for possible relaxations of these hypotheses. 

\begin{proposition}
Assume hypotheses (H1)--(H4), with committee size $n$, and competence $\theta$.
Then the probability that the votes result in a particular table $T=(x,y,z,t)$ is, under
the different states of nature, 
\begin{align}
\label{eq:PandQ}
\text{if $P\wedge Q$, then} \quad& \binom{n}{x,y,z,t}\theta^{2x+y+z}(1-\theta)^{y+z+2t}\ ,
\\ 
\label{eq:PandnoQ}
\text{if $P\wedge \neg Q$, then} \quad& \binom{n}{x,y,z,t}\theta^{x+2y+t}(1-\theta)^{x+2z+t}\ ,
\\ 
\label{eq:noPandQ}
\text{if $\neg P\wedge Q$, then} \quad& \binom{n}{x,y,z,t}\theta^{x+2z+t}(1-\theta)^{x+2y+t}\ ,
\\ 
\label{eq:noPandnoQ}
\text{if $\neg P\wedge \neg Q$, then} \quad& \binom{n}{x,y,z,t}\theta^{y+z+2t}(1-\theta)^{2x+y+z}\ ,
\end{align}
where $!$ means the factorial of a number. 
\end{proposition}
\begin{proof}
See Proposition A.4 of \cite{AlabertFarre2022}.
\end{proof}
Let us denote $\mathbb P_{P\wedge Q}$ the probabilities computed under the state of nature
$P\wedge Q$. According to the proposition above, the probability of obtaining a true
positive when rule $r$ is employed is the sum of the probabilities (\ref{eq:PandQ}) for 
all tables $T$ such that $r(T)=1$:
\begin{equation*}\label{probaTP}
\mathbb P_{P\wedge Q}(\text{TP}) 
  = \sum_{\{r(x,y,z,t)=1\}} \binom{n}{x,y,z,t}\theta^{2x+y+z}(1-\theta)^{y+z+2t} \ .
\end{equation*}
Therefore, the probability of incurring a false negative is 
\begin{equation*}\label{probaFN}
\mathbb P_{P\wedge Q}(\text{FN}) 
  = 1-\mathbb P_{P\wedge Q}(\text{TP})
  = \sum_{\{r(x,y,z,t)=0\}} \binom{n}{x,y,z,t}\theta^{2x+y+z}(1-\theta)^{y+z+2t} 
  \ .
\end{equation*}

We cannot proceed in a completely analogous way to define true negatives and false positives,
because $\neg(P\wedge Q)$ is not a state of nature, but an ensemble of three states, each
of which may yield different probabilities. At this point, there are two possible modelling
paths, according to the information available: Either there is no further information about 
the true state of nature (or we do not want to use it); or, there is enough information to
postulate an ``a priori'' probability $\pi$ on the states of nature, and we can follow a Bayesian approach. 

The main line
in this paper is the first path, always applicable. Let us deviate for a moment and 
sketch the second one, which corresponds to a situation considered, among others, 
in Terzopoulou and Endriss \cite{TerzopoulouEndrissSAGT2019},
Bovens and Rabinowicz \cite{Bovens2006} and Bozbay \cite{Bozbay2019}:
In the Bayesian approach, $\mathbb P_{P\wedge Q}$ is interpreted as a conditional 
probability given $P\wedge Q$, and analogously for the other ones, that we denote 
$\mathbb P_{P\wedge \neg Q}$, $\mathbb P_{\neg P\wedge Q}$, and
$\mathbb P_{\neg P\wedge\neg Q}$. Hence, the probability of a true negative
in this setting will be
\begin{equation}\label{probaTN}
\begin{split}
\mathbb P_{\neg(P\wedge Q)}(\text{TN})= &
\mathbb P_{P\wedge\neg Q}(\text{TN})\cdot \pi(P\wedge\neg Q)
\\ &+
\mathbb P_{\neg P\wedge Q}(\text{TN})\cdot \pi(\neg P\wedge Q)
\\ &+
\mathbb P_{\neg P\wedge \neg Q}(\text{TN})\cdot \pi(\neg P\wedge \neg Q)
\ ,
\end{split}
\end{equation}
and then the probability of a false positive is given by
\begin{equation}\label{probaFP}
  \mathbb P_{\neg(P\wedge Q)} (\text{FP}) = 1- \mathbb P_{\neg(P\wedge Q)} (\text{TN})
\end{equation}
For example, if $\pi$ is assumed to give the same probability to all three
negative states, then $\mathbb P_{\neg(P\wedge Q)}(\text{TN})$ will be the arithmetic mean
of the probabilities of TN under each state. This is the chosen prior in \cite{TerzopoulouEndrissSAGT2019};
that of \cite{Bovens2006} is different, and Bozbay \cite{Bozbay2019} completely forbids
the result $\neg P \wedge \neg Q$. In general, if the committee knows the prior, the 
independence in the judgments of the premisses (H4) cannot be assumed.

Now, using expressions 
(\ref{eq:PandnoQ}), (\ref{eq:noPandQ}) and (\ref{eq:noPandnoQ}), the probabilities of a true
negative under
the three negative states are
\begin{align}
\label{probaTN1}
\mathbb P_{P\wedge \neg Q}(\text{TN}) 
  = \sum_{\{r(x,y,z,t)=0\}} \binom{n}{x,y,z,t}\theta^{x+2y+t}(1-\theta)^{x+2z+t}
\\ 
\label{probaTN2}
\mathbb P_{\neg P\wedge Q}(\text{TN}) 
  = \sum_{\{r(x,y,z,t)=0\}} \binom{n}{x,y,z,t}\theta^{x+2z+t}(1-\theta)^{x+2y+t}
\\  
\label{probaTN3}
\mathbb P_{\neg P\wedge \neg Q}(\text{TN}) 
  = \sum_{\{r(x,y,z,t)=0\}} \binom{n}{x,y,z,t}\theta^{y+z+2t}(1-\theta)^{2x+y+z} 
\end{align}
and the probability of a false positive is then computed from (\ref{probaTN}) and (\ref{probaFP}).

After this digression, let us turn to our main setting. 
For the non-Bayesian situation, we can resort to the following analogy with the classical 
theory of Statistical Hypothesis Testing:
Suppose one has to decide if there is enough evidence that a certain population parameter
is equal to a value $C$, as provided by a sample drawn from the population. 
To this end, one computes how likely
the observed sample could have been produced by the value of the parameter  
in the complement set $\neg C$ which is ``the closest'' to $C$.
If that likelihood is acceptable (by some numerical threshold), the decision is to stick
to the ``null'' (status quo) conclusion $\neg C$. If it is not acceptable, $C$ is proclaimed as
the new estimated conclusion.

Translating the analogy to our case, we must ask ourselves which of the states of nature
belonging to the complement set $\neg(P\wedge Q)$ is closest to $P\wedge Q$. Intuitively, 
$P\wedge\neg Q$ and $\neg P\wedge Q$ are equally close, and are closer than
$\neg P\wedge\neg Q$.
This is rigorously stated in the next proposition.
Although intuitive, the rigorous proof is a 
little bit technical. We use a probabilistic procedure called \emph{coupling}, that
transforms inequalities about probabilities into inequalities about random variables. 

These arguments support the \emph{definition} of the ``probability'' of a false positive 
as the probability that a rule decides $C$ when $P\wedge\neg Q$ is the true state of nature: 
$\mathbb P_{P\wedge\neg Q}(\text{FP}):=1-\mathbb P_{P\wedge\neg Q}(\text{TN})$.
It can be also thought as the maximum of the probabilities of a false positive for all
possible choices of a prior $\pi$ on the states of nature. It is therefore a conservative
estimate of the possible error, in response to the lack of information about the underlying truth.

\begin{proposition}
Under hypothesis (H1)-(H4), for any admissible decision rule $r$, 
\begin{equation*}
\mathbb P_{P\wedge \neg Q} \{r=1\} =
\mathbb P_{\neg P\wedge Q} \{r=1\} \ge \mathbb P_{\neg P \wedge \neg Q} \{r=1\}
\ .
\end{equation*}
\end{proposition} 
\begin{proof}
The first equality is clear from Definition \ref{def:admissible}, item 1. 
We only need to prove the inequality on the right.
Let ${\mathbb M}_{\neg P\wedge\neg Q}$ and ${\mathbb M}_{P\wedge\neg Q}$
be two probability measures defined on the subsets of a set $\Omega$, and let $T$ be a random
variable $T\colon\Omega\rightarrow{\mathbb T}$ such that the law of $T$ under 
${\mathbb M}_{\neg P\wedge\neg Q}$ is ${\mathbb P}_{\neg P\wedge\neg Q}$ 
and the law of $T$ under ${\mathbb M}_{P\wedge\neg Q}$ is ${\mathbb P}_{P\wedge\neg Q}$.
That is, using (\ref{eq:noPandnoQ}) and (\ref{eq:PandnoQ}), 
\begin{align*}
{\mathbb M}_{\neg P\wedge\neg Q}
\{\omega:\ T(\omega)=(x,y,z,t)\}
&=
{\mathbb P}_{\neg P\wedge\neg Q}\{(x,y,z,t)\}= 
\binom{n}{x,y,z,t}\theta^{y+z+2t}(1-\theta)^{2x+y+z}\ ,
\\
{\mathbb M}_{P\wedge\neg Q}
\{\omega:\ T(\omega)=(x,y,z,t)\}
&=
{\mathbb P}_{P\wedge\neg Q}\{(x,y,z,t)\}= 
\binom{n}{x,y,z,t}\theta^{x+2y+t}(1-\theta)^{x+2z+t}\ .
\end{align*}

Suppose we could define another random variable $S\colon\Omega\rightarrow \mathbb T$ such that 
\begin{enumerate}[a)]
\item
$T(\omega)\le S(\omega)$, for all $\omega\in\Omega$, and
\item 
The law of $S$ under ${\mathbb M}_{\neg P\wedge\neg Q}$ 
coincides with the law of $T$ under ${\mathbb M}_{P\wedge\neg Q}$.
\end{enumerate}
Then, since  $r$ is order-preserving, we will have
$
\{\omega:\ r(S(\omega))=1\}
\supseteq
\{\omega:\ r(T(\omega))=1\}
$, 
 and the conclusion
\begin{equation*}
{\mathbb P}_{P\wedge\neg Q}\{r(T)=1\}=
{\mathbb P}_{\neg P\wedge\neg Q}\{r(S)=1\}\ge
{\mathbb P}_{\neg P\wedge\neg Q}\{r(T)=1\}
\ .
\end{equation*}
Let us prove the existence of $S\colon\Omega\rightarrow \mathbb T$ with the properties
a) and b) above, and we are done: Let $T_1,\dots, T_n$ be independent identically distributed 
random variables $T_i\colon\Omega\rightarrow\mathbb T$ with the same law as $T$ but 
for the vote of one individual. 
We will switch to table notation again, for clarity, in the 
rest of the proof.

Let $S_1,\dots, S_n$ be another collection of independent identically distributed random
variables $S_i\colon\Omega\rightarrow\mathbb T$, defined as follows:
\begin{align*}
&\text{If $T_i=\Big[\,\votes{1}{0}{0}{0}\,\Big]$, 
   then $S_i=\Big[\,\votes{1}{0}{0}{0}\,\Big]$}  
\\   
&\text{If $T_i=\Big[\,\votes{0}{1}{0}{0}\,\Big]$, 
   then $S_i=\Big[\,\votes{0}{1}{0}{0}\,\Big]$} 
\\    
&\text{If $T_i=\Big[\,\votes{0}{0}{1}{0}\,\Big]$, 
   then $S_i
   =\begin{cases} 
     \Big[\,\votes{0}{0}{1}{0}\,\Big]
     \text{ with probability } \tfrac{1-\theta}{\theta}
     \\[5pt] 
     \Big[\,\votes{1}{0}{0}{0}\,\Big]
     \text{ with probability } \tfrac{2\theta-1}{\theta}
    \end{cases}
    $}
\\    
&\text{If $T_i=\Big[\,\votes{0}{0}{0}{1}\,\Big]$, 
   then $S_i
   =\begin{cases} 
     \Big[\,\votes{0}{0}{0}{1}\,\Big]
     \text{ with probability } \tfrac{1-\theta}{\theta}
     \\[5pt] 
     \Big[\,\votes{0}{1}{0}{0}\,\Big]
     \text{ with probability } \tfrac{2\theta-1}{\theta}
    \end{cases}
    $}    
\end{align*}
We have clearly that $T_i(\omega)\le S_i(\omega)$, for all $\omega\in\Omega$.

Let us compute the law of $S_i$ under ${\mathbb M}_{\neg P\wedge\neg Q}$, using conditional
probabilities to the value of $T_i$:
\begin{align*}
{\mathbb M}_{\neg P\wedge\neg Q}
\Big\{S_i=\Big[\,\votes{1}{0}{0}{0}\,\Big]\Big\}
= &
{\mathbb M}_{\neg P\wedge\neg Q}
\Big\{
  \condprob{S_i=\Big[\,\votes{1}{0}{0}{0}\,\Big]}
           {T_i=\Big[\,\votes{1}{0}{0}{0}\,\Big]}
\Big\}
\cdot
{\mathbb M}_{\neg P\wedge\neg Q}
\Big\{
  T_i=\Big[\,\votes{1}{0}{0}{0}\,\Big]
\Big\}
\\
&+
{\mathbb M}_{\neg P\wedge\neg Q}
\Big\{
  \condprob{S_i=\Big[\,\votes{1}{0}{0}{0}\,\Big]}
           {T_i=\Big[\,\votes{0}{0}{1}{0}\,\Big]}
\Big\}
\cdot
{\mathbb M}_{\neg P\wedge\neg Q}
\Big\{
  T_i=\Big[\,\votes{0}{0}{1}{0}\,\Big]
\Big\}
\\
=&
(1-\theta)^2+\frac{2\theta-1}{\theta}\cdot\theta(1-\theta)=\theta(1-\theta)
\ .
\\
{\mathbb M}_{\neg P\wedge\neg Q}
\Big\{S_i=\Big[\,\votes{0}{1}{0}{0}\,\Big]\Big\}
= &
{\mathbb M}_{\neg P\wedge\neg Q}
\Big\{
  \condprob{S_i=\Big[\,\votes{0}{1}{0}{0}\,\Big]}
           {T_i=\Big[\,\votes{0}{1}{0}{0}\,\Big]}
\Big\}
\cdot
{\mathbb M}_{\neg P\wedge\neg Q}
\Big\{
  T_i=\Big[\,\votes{0}{1}{0}{0}\,\Big]
\Big\}
\\
&+
{\mathbb M}_{\neg P\wedge\neg Q}
\Big\{
  \condprob{S_i=\Big[\,\votes{0}{1}{0}{0}\,\Big]}
           {T_i=\Big[\,\votes{0}{0}{0}{1}\,\Big]}
\Big\}
\cdot
{\mathbb M}_{\neg P\wedge\neg Q}
\Big\{
  T_i=\Big[\,\votes{0}{0}{0}{1}\,\Big]
\Big\}
\\
=&
\theta(1-\theta) + \frac{2\theta-1}{\theta}\cdot\theta^2= \theta^2
\ .
\\
{\mathbb M}_{\neg P\wedge\neg Q}
\Big\{S_i=\Big[\,\votes{0}{0}{1}{0}\,\Big]\Big\}
= &
{\mathbb M}_{\neg P\wedge\neg Q}
\Big\{
  \condprob{S_i=\Big[\,\votes{0}{0}{1}{0}\,\Big]}
           {T_i=\Big[\,\votes{0}{0}{1}{0}\,\Big]}
\Big\}
\cdot
{\mathbb M}_{\neg P\wedge\neg Q}
\Big\{
  T_i=\Big[\,\votes{0}{0}{1}{0}\,\Big]
\Big\}
\\
=&
\frac{1-\theta}{\theta}\cdot \theta(1-\theta) = (1-\theta)^2
\ .
\\
{\mathbb M}_{\neg P\wedge\neg Q}
\Big\{S_i=\Big[\,\votes{0}{0}{0}{1}\,\Big]\Big\}
= &
{\mathbb M}_{\neg P\wedge\neg Q}
\Big\{
  \condprob{S_i=\Big[\,\votes{0}{0}{0}{1}\,\Big]}
           {T_i=\Big[\,\votes{0}{0}{0}{1}\,\Big]}
\Big\}
\cdot
{\mathbb M}_{\neg P\wedge\neg Q}
\Big\{
  T_i=\Big[\,\votes{0}{0}{0}{1}\,\Big]
\Big\}
\\
=&
\frac{1-\theta}{\theta}\cdot \theta^2 = \theta(1-\theta) 
\ .
\end{align*} 

We see that the law of each $S_i$ under $\mathbb M_{\neg P\wedge\neg Q}$ coincides
with that of $T_i$ under  $\mathbb M_{P\wedge\neg Q}$. Now, we have that
\begin{equation*}
T(\omega):=\sum_{i=1}^n T_i(\omega)\le \sum_{i=1}^n S_i(\omega)=:S(\omega)\ ,\quad \omega\in\Omega
\ ,
\end{equation*}
and that 
$T$ has the law given by (\ref{eq:noPandnoQ}) under $\mathbb M_{\neg P\wedge\neg Q}$,
and by (\ref{eq:PandnoQ}) under $\mathbb M_{P\wedge\neg Q}$,
and $S$ has the law given by (\ref{eq:PandnoQ}) under $\mathbb M_{\neg P\wedge\neg Q}$.

Hence, $T$ and $S$ are the random variables we were looking for, and the proof is complete.
\end{proof}
As a corollary, since the maximum 
is always greater than any weighted mean in (\ref{probaTN}), we get that $\mathbb P_{P\wedge\neg Q}\{r=1\}$
is greater or equal than the probability of a false positive computed with any prior distribution on the
set $\neg(P\wedge Q)$.

In a completely analogous way, one can also prove that for admissible rules 
the probability to conclude $C$ is greater when the true state of nature is $P\wedge Q$ than
with any other state. 

We do not need any more the subindexes in the probabilities of false positives and false 
negatives, since $\mathbb P(\text{FN})$ always refers to the state $P\wedge Q$, 
and $\mathbb P(\text{FP})$ refers to the state $P\wedge\neg Q$ (or to the 
given prior on $\neg(P\wedge Q)$ in the Bayesian case). 
Instead, we will subindex $\mathbb P$ by the rule employed. For reference
in the sequel, we repeat here the formulae for FP and FN: For any rule $r\colon {\mathbb T}\rightarrow \{0,1\}$, 
\begin{align}
\label{eq:PFP}
\mathbb P_r(\text{FP})
  &= 
  \sum_{\{r(x,y,z,t)=1\}} \binom{n}{x,y,z,t}\theta^{x+2y+t}(1-\theta)^{x+2z+t} \ ,  
  \\  
\label{eq:PFN}
\mathbb P_r(\text{FN}) 
  &= \sum_{\{r(x,y,z,t)=0\}} \binom{n}{x,y,z,t}\theta^{2x+y+z}(1-\theta)^{y+z+2t} 
  \ .  
\end{align}

\subsection{Optimisation criterion}\label{ssec:OptCrit}
We want to obtain the best decision rule, under the probabilistic model stated above 
and the optimisation criterion what we develop in this subsection. This criterion was
introduced in \cite{AlabertFarre2022} and is based on minimising a weighted sum
of the probability to commit a false positive and the probability to commit a 
false negative. It can be thought as a multi-objective optimisation problem, but that
point of view does not contribute any practical insight.  

Any rule $r\colon \mathbb T\rightarrow\{0,1\}$ (admissible or not) has associated probabilities 
of producing a False Positive ${\mathbb P}_r(\text{FP})$ and a False Negative ${\mathbb P}_r(\text{FN})$
according to formulae (\ref{eq:PFP})--(\ref{eq:PFN}). Despite the simplified notation, recall that
these two probabilities stem from different states of nature. If both failures are
considered equally harmful, it is natural to look for the admissible rule $r\in{\cal A}$ that
minimises the sum   
\begin{equation}\label{eq:AOT}
{\mathbb P}_r(\text{FP}) + {\mathbb P}_r(\text{FN})\ .
\end{equation}
If one of them is considered worse that the other, one can take a weighted sum
\begin{equation}\label{eq:WAOT}
L_w(r):=w{\mathbb P}_r(\text{FP}) + (1-w) {\mathbb P}_r(\text{FN})\ ,
\end{equation}
where $w$ is a real number, $0<w<1$, as the loss function that to minimise. For example, 
if a false positive is deemed 
twice as harmful as a false negative, $w=\frac{2}{3}$ is the suitable value. 
Note that the weight $w$ is a modelling choice relative to each particular application, 
and it is supposed to be fixed in advance of the voting stage.

In Statistics, the combination (\ref{eq:AOT}) of probabilities of the two types of error 
is called the \emph{area of the triangle},
a term that comes from its origin in signal processing and the graphical 
methodology called \emph{Receiving
Operating Characteristics} (ROC). We refer the reader to Fawcett 
\cite{Fawcett:2006:IRA:1159473.1159475} for a simple introduction to ROC.

If $r$ and $s$ are two admissible rules, and $r\le s$ 
(equivalently, $U_r\subseteq U_s$, where $U_r$ and $U_s$
are the upper sets defining $r$ and $s$ respectively), then obviously, from (\ref{eq:PFP})
and (\ref{eq:PFN}),
\begin{align*}
\mathbb P_r(\text{FP}) &\le \mathbb P_s(\text{FP}) \\
\mathbb P_r(\text{FN}) &\ge \mathbb P_s(\text{FN}) \ .
\end{align*}
This means that $r\mapsto \mathbb P_r(\text{FP})$ and $r\mapsto \mathbb P_r(\text{FN})$ 
are respectively an increasing function and a decreasing function 
defined on the poset of admissible rules $(\cal A,\le)$.
Moreover the rule 
$r\equiv 0$ (always conclude $\neg C$)
satisfies $\mathbb P_r(\text{FP})=0$ and $\mathbb P_r(\text{FN})=1$, and the rule 
$r\equiv 1$ (always conclude $C$)
satisfies $\mathbb P_r(\text{FP})=1$ and $\mathbb P_r(\text{FN})=0$.

As we said before, the value $\theta=\frac{1}{2}$ can be excluded because it is trivial: 
If $w<\frac{1}{2}$,
the decision must be always $C$; if $w>1/2$, the decision must be $\neg C$; and 
if $w=1/2$, then the problem is completely equivalent to a single coin toss. 
See the appendix for the details.

\section{Main results}\label{sec:Main}

Suppose we have an admissible rule $r$, with sets $\mathbb T^+$ of positive tables 
and $\mathbb T^0$ of null tables (recall the definitions in Section 2). 
If we choose a table $T\in \mathbb T^0$ and move it to $\mathbb T^+$, we are 
increasing the probability of a false positive and at the same time decreasing the
probability of a false negative. In doing that, we also have to move
its transposed table, to maintain the first condition of admissibility. 
This movement may result in a decrease or an increase of the loss function $L_w$. 

\begin{definition}
Let $r\colon {\mathbb T}\rightarrow \{0,1\}$ be an admissible decision rule, 
with positive set $\mathbb T^+$ and null set $\mathbb T^0$. 
If moving a table $T$ and its transposed table (and no other) from $\mathbb T^0$ to $\mathbb T^+$ 
results in 
a decrease of the loss function, we will say that we have a \emph{good table}.  
\end{definition}

Here ``good'' only means that the voting table $T$ ``should be supporting decision C''.
Thus, it seems that the set $\mathbb T^+$ of the optimal rule must consist of the good tables 
and no others.
However, we still have to see that the rule defined in this way 
is indeed admissible. 

\begin{theorem}
The rule whose positive set $\mathbb T^+$ consists exactly of the good tables 
is admissible and optimal.
\end{theorem} 
This is a consequence of the following two lemmas, which are interesting in their own.
Lemma \ref{lem:goodTable} characterises the good tables in terms of $w$ and $\theta$,
and confirms that a table and its transposed are both good or both bad. 
Lemma \ref{lem:goodUpper} proves that the second condition of admissibility is
also satisfied, in view of Proposition \ref{pr:goodUpper}.

\begin{lemma}\label{lem:goodTable}
Given weight $0<w<1$ and competence $\frac{1}{2}<\theta<1$, 
a table $T=(x,y,z,t)\in \mathbb T$ is a good table if and only if 
\begin{equation}\label{eq:goodTable}
\Big(\frac{\theta}{1-\theta}\Big)^{(y-z)-(x-t)}
+
\Big(\frac{\theta}{1-\theta}\Big)^{(z-y)-(x-t)}
< \frac{2(1-w)}{w}\ . 
\end{equation}  
\end{lemma}
\begin{proof}
Using formulae (\ref{eq:PFP}) and (\ref{eq:PFN}), 
the change in the quantity $w\mathbb P(\text{FP})+(1-w)\mathbb P(\text{FN})$
when $T$ and its transposed table are moved from the null to the positive set is given by
\begin{equation}\label{FP-FN}
\begin{split}
\binom{n}{x,y,z,t}\Big[
&
w\theta^{x+2y+t}(1-\theta)^{x+2z+t}+w\theta^{x+2z+t}(1-\theta)^{x+2y+t}
\\ 
&
- 2(1-w)\theta^{2x+y+z}(1-\theta)^{y+z+2t}\Big]
\end{split}
\end{equation}
Dropping the factorials and dividing by 
\begin{equation*}
\theta^{2x+y+z}(1-\theta)^{y+z+2t} 
\ ,
\end{equation*} we see
that (\ref{FP-FN}) is negative when
\begin{equation*}
\theta^{y+t-x-z}(1-\theta)^{x+z-y-t}+\theta^{z+t-x-y}(1-\theta)^{x+y-z-t} < \frac{2(1-w)}{w}
\ ,
\end{equation*}
and this is immediately equivalent to (\ref{eq:goodTable}). In the case $y=z$, the quantity
(\ref{FP-FN}) should be divided by two, but the conclusion is the same.
\end{proof}
\begin{lemma}\label{lem:goodUpper}
The set of good tables is an upper set of $(\mathbb{T},\le)$. 
\end{lemma}
\begin{proof}
It is easy to see that, for each fixed $\frac{1}{2}<\theta<1$, the function 
\begin{equation}\label{eq:Hfunction}
H(x,y,t,z)= 
\Big(\frac{\theta}{1-\theta}\Big)^{(y-z)-(x-t)}
+
\Big(\frac{\theta}{1-\theta}\Big)^{(z-y)-(x-t)}
\end{equation}
is decreasing in $(\mathbb T,\le)$. It is enough to check it for the pairs of the transitive reduction
(\ref{eq:transreduc}). Moving from the smallest to the greatest table of the pair, one of the terms 
in (\ref{eq:Hfunction}) remains unchanged whereas 
the other one is multiplied by $\big(\frac{\theta}{1-\theta}\big)^{-2}<1$. Therefore, if a table is good, according to 
Lemma \ref{lem:goodTable} a greater table in the poset $(\mathbb T,\le)$ is also good, and the good
tables indeed form an upper set. 
\end{proof}

Notice that the optimal rule according to $L_w$ among those satisfying the first 
admissibility condition only, automatically satisfies the second. 

\bigskip
The next theorem determines when, under the optimal decision rule, a voting table leads to a verdict of
$C=P\wedge Q$ or the opposite, for the symmetric case $w=1/2$. This is a further characterisation of
the condition (\ref{eq:goodTable}). We prove this special case because
the statement and the proof are neater, and the extension to general $0<w<1$ is straightforward, as will be seen
after the theorem.

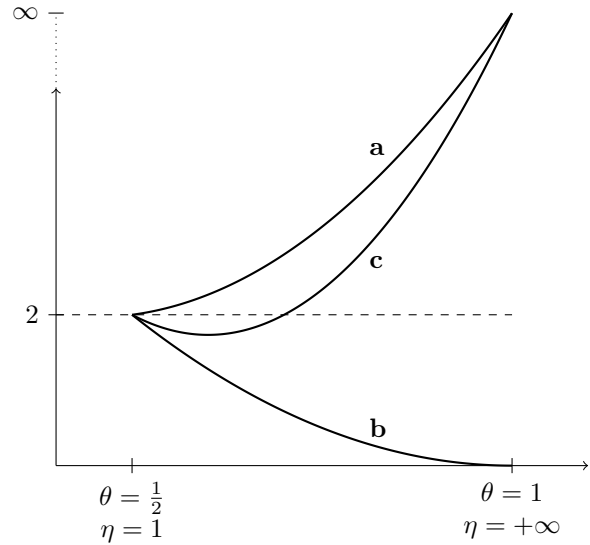
\begin{figure}[h]
{\centering\small
\begin{tikzpicture}[scale=1]  
\begin{scope}[shift={(1,1)}]
  \draw[->] (0,0) -- (7,0);
  \draw[->] (0,0) -- (0,5);
  \draw[-, dotted] (0,5) -- (0,6);
  \draw[-] (-0.1,6) -- (0.1,6); 
  \node[left] at (-0.1,6) {$\infty$};
  \draw[-] (-0.1,2) -- (0.1,2); 
  \node[left] at (-0.1,2) {$2$};
  \draw[-] (1,-0.1) -- (1,0.1); 
  \node[below, align=center] at (1,-0.1) {$\theta=\frac{1}{2}$\\ $\eta=1$};
  \draw[-] (6,-0.1) -- (6,0.1); 
  \node[below, align=center] at (6,-0.1) {$\theta=1$\\ $\eta=+\infty$};
  \draw[-, dashed] (0,2) -- (6,2);
  \node[above right] at (4,4) {\textbf{a}};
  \node[above right] at (4,0.25) {\textbf{b}};
  \node[above right] at (4,2.5) {\textbf{c}};
  \begin{scope}[shift={(1,2)}]
    \draw[thick, x=1cm, y=1cm] 
    plot[domain=0:5, samples=144, smooth] (\x, {(2/15)*\x^2+(2/15)*\x});
    \draw[thick, x=1cm, y=1cm] 
    plot[domain=0:5, samples=144, smooth] (\x, {(2/25)*\x^2-(20/25)*\x});
    \draw[thick, x=1cm, y=1cm] 
    plot[domain=0:5, samples=144, smooth] (\x, {(4/15)*\x^2-(8/15)*\x});
  \end{scope}
\end{scope}
\end{tikzpicture}
\caption{Different possible shapes for function $G$ of (\ref{eq:G}) according to the 
three cases in the proof of Theorem \ref{thm:optrule}.}
\label{fig:Gfunctions}
}
\end{figure}

In words, the theorem says that: if votes in favour of $P\wedge Q$
are less than those in favour of $\neg P\wedge\neg Q$ or there is a tie, then the decision 
must be $\neg C$; otherwise, if the difference is greater than the difference in absolute value between 
votes for $P\wedge \neg Q$ and votes for $\neg P\wedge Q$, then the decision must be $C$; 
otherwise, the decision must be $C$ if the competence $\theta$ of the committee is below a
certain threshold (which can be computed to any desired accuracy), and $\neg C$ if it is above.
\begin{theorem}\label{thm:optrule}
Assume $w=1/2$, and let $\theta$ be any competence level, $\frac{1}{2}<\theta<1$. Assume $r$ is an optimal rule
and denote 
\begin{equation*}
\rho:=x-t \text{ \ and \ } \alpha:=|y-z| \ .
\end{equation*} 
Then, given a voting table $T=(x,y,z,t)$, 
\begin{enumerate}[a)]
\item 
  if $\rho \le 0$, then $r(T)=0$.
\item 
  if $\rho>\alpha$, then $r(T)=1$.
\item 
  if $0< \rho < \alpha$, then there exists $\theta_0\in(\frac{1}{2},1)$ such that 
  $r(T)=1$ for $\theta<\theta_0$, and $r(T)=0$ for $\theta>\theta_0$.
\end{enumerate}
 And these are all possible cases. 
\end{theorem}
\begin{proof}
Starting with the last claim, these are all possible cases because 
$n=x+y+z+t$ odd implies that $\alpha$ and $\rho$ have different parity. In 
particular, $\rho\neq\alpha$ and $\rho\neq -\alpha$.
 
Consider the bijective increasing transformation  
$\eta=\frac{\theta}{1-\theta}$, which maps 
$(\frac{1}{2},1)$ onto $(1,\infty)$. 
The left-hand side of (\ref{eq:goodTable}) can then be written as a function 
\begin{equation}\label{eq:G}
G(\eta):=\eta^{-\rho-\alpha}+\eta^{-\rho+\alpha}\ ,\quad
\eta\in(1,\infty)\ ,
\end{equation}
with $\rho$ and $\alpha$ integers, $\alpha\ge 0$. The restriction
$x+y+z+t=n$ implies that 
\begin{equation*}
\begin{aligned}
          0 & \le\alpha\le n \\
-(n-\alpha) & \le\rho\le n-\alpha 
\end{aligned}
\end{equation*}
with $\rho$ and $\alpha$ of different parity, as already noted.
According to Lemma \ref{lem:goodTable}, the table $T$ is good for values of $\eta$ such that $G(\eta)<2$,
and the optimal rule $r$ should assign it the value 1; for values of $\eta$ such that $G(\eta)>2$, the table
is bad and we must have $r(T)=0$.

Clearly, $G$ is differentiable on $(1,\infty)$ and $\lim_{\eta\searrow 1} G(\eta) = 2$, 
for all $\rho$ and $\alpha$.
Also, the derivative of $G$ can always be written as
\begin{equation}\label{eq:G'}
G'(\eta)=\eta^{-\rho-1}\Big[(\alpha-\rho)\eta^{\alpha}-(\alpha+\rho)\eta^{-\alpha}\Big]\ ,
\end{equation}
and we have $\lim_{\eta\downarrow 1} G'(\eta)=-2\rho$. Thus, $G$ takes the value 2
at the left boundary of
the domain of interest, and starts from there decreasing or increasing according to the 
sign of $\rho$. 

Let us now proceed with the three cases of the statement. Please refer to Figure 
\ref{fig:Gfunctions}.
\begin{description}
\item [Case a:] $\rho\le 0$. \\
If $\rho=0$, then $\alpha$ cannot be zero, by parity. We have $\alpha\ge 1$ and $G$ is
clearly increasing, with $\lim_{\eta\to\infty} G(\eta)=+\infty$. If $\rho<0$, $G$ is also
increasing for all $\alpha\ge 0$ and $\lim_{\eta\to\infty} G(\eta)=+\infty$ again. The table
is bad in both situations.

\item [Case b:] $\rho>\alpha$. \\
Both exponents in (\ref{eq:G}) are negative, and $G$ is therefore decreasing. 
Moreover $\lim_{\eta\to\infty} G(\eta)=0$. The table is good.

\item [Case c:] $0<\rho<\alpha$. \\
Now $G'$ is negative near $\eta=1$; therefore $G$ is decreasing at least on some
interval to the right of 1. Solving for $\eta$ in $G'(\eta)=0$, and taking into 
account that $\alpha+\rho>\alpha-\rho>0$, we find that 
\begin{equation}\label{eq:minG}
\eta^* = \Big(\frac{\alpha+\rho}{\alpha-\rho}\Big)^{1/2\alpha}\in (1,\infty)
\end{equation}
is the only critical point of $G$, which is a minimum since $\lim_{\eta\to\infty} G(\eta)=+\infty$.

Thus, there exists a unique point $\eta_0\in (\eta^*,\infty)$ where $G(\eta_0)=2$. 
The table is good for $\eta\in (1,\eta_0)$, and bad for $\eta\in (\eta_0,\infty)$.  

In other words, for a competence value $\theta$ less than $\theta_0:=\frac{\eta_0}{1+\eta_0}$,
the table is good. For $\theta>\theta_0$, the table is bad. This finishes the proof.
\end{description}
\end{proof}

The case for general $0<w<1$ is very easy to explain with the help of Figure \ref{fig:Gfunctions}.
For $w<\frac{1}{2}$, the dashed horizontal line is above level 2. All tables are good for
small enough competence levels $\theta$ (the decision must always be $C$). 
The tables of type \textbf{a} and \textbf{c} are bad for 
$\theta$ greater than some $\theta_0$. The tables of type \textbf{b} are all good for all 
competence levels.

For $w>\frac{1}{2}$, the dashed line is below level 2, and all tables are bad for low 
enough competence levels (the decision must always be $\neg C$). Tables \textbf{a} 
stay bad for all $\theta$, and tables \textbf{b} turn good after some point. For tables of type \textbf{c}, 
two things may happen: Either they are always bad or, as $\theta$ increases, they have 
an interval of ``goodness'' before turning bad again. 

All intersection points are easily computed to any desired precision by solving numerically
for $\eta$ the equation $G(\eta)=\frac{2(1-w)}{w}$ on $(1,\infty)$. See Section \ref{sec:ComputSoft}
for more details.

Lemma \ref{lem:goodTable} allows a notable conceptual, notational and computational
simplification: 
Since function $G$ in (\ref{eq:G}) only depends on $\rho=x-t$ and $\alpha=|y-z|$,
the tables in $\mathbb T$ with the same $\rho$ and $\alpha$ will all be good or bad, 
once $\theta$ and $w$ are fixed. 
If, on the contrary, two given tables do not
share these values, they produce two different functions $G$.

This allows to consider an equivalence relation in $(\mathbb T,\le)$ that gives rise to 
a quotient ranked poset, reducing in this way the complexity of the Hasse diagram and 
the computations.
Define
\begin{equation*}
(x,y,z,t)\sim(x',y',z',t') \Longleftrightarrow x-t=x'-t'\text{ and } |y-z|=|y'-z'| 
\ .
\end{equation*}
In particular, this equivalence relation identify transposed tables. 
The elements of the quotient set $\mathbb T /\negthickspace\sim$ are classes of voting tables 
and can be represented by the pair $(\rho, \alpha)$. We can write $T\in(\rho,\alpha)$
if $T$ is in the class represented by $(\rho,\alpha)$.

Now define the preorder relation in $\mathbb T /\negthickspace\sim$ given by 
\begin{align*}
(\rho,\alpha)\le(\rho',\alpha') \Longleftrightarrow 
&\text{ there exist $T\in(\rho,\alpha)$ and $T'\in(\rho',\alpha')$ such that $T\le T'$} 
\ .
\end{align*}
We use the same symbol '$\le$' for both relations in $\mathbb T$ and 
$\mathbb T/\negthickspace\sim$, 
since there is no possible confusion. It can be proved in general that a relation
defined in this way in the quotient set is reflexive and transitive, therefore
a \emph{preorder}. In general it is not antisymmetric. 

We shall prove that in our case the antisymmetry holds, so that we have again
a partial order. To this end, we make use of the following lemma 
(see Hallam \cite{Hallam2015}).
A proof is included in the appendix, for the reader convenience.
 
\begin{lemma}\label{lem:preorder2poset}
Let $(X,\le)$ be a finite poset, $\sim$ an equivalence relation on $X$, and 
the preorder $\preceq$ on $X/\negthickspace\sim$ defined as: 
$\bar x\preceq\bar y \Leftrightarrow \exists x\in\bar x, \exists y\in\bar y:\ x\le y$.\\
Assume that if $\bar x\preceq\bar y$ in $X/\negthickspace\sim$, then for all 
$x\in \bar x$, there exists
$y\in \bar y$ such that $x\le y$ in $X$.
Then, $(X/\negthickspace\sim, \preceq)$ is a poset.    
\end{lemma}
It is not difficult to show that the hypothesis of the lemma holds in our case; 
see the appendix.

Since we are identifying tables in the same rank level, 
the resulting quotient poset is also ranked, with the same rank function  $\rho$.

\begin{figure}[t]
{\footnotesize
\begin{minipage}[b][0.6\textheight][c]
{0.5\textwidth}  
{\centering
\begin{tikzpicture}[scale=0.8]
    {\color{black!50!white}
    \node (A) at (0,0) {$(-3,0)$};
    \node (B) at (0,1) {$(-2,1)$};
    \node (C) at (-1,2) {$(-1,0)$};
    \node (D) at (1,2) {$(-1,2)$};
    \node (E) at (-1,3) {$(0,1)$};
    \node (F) at (1,3) {$(0,3)$}; }
    \node (G) at (-1,4) {$\mathbf{(1,0)}$};
    \node (H) at (1,4) {$(1,2)$};
    \node (I) at (0,5) {$\mathbf{(2,1)}$};
    \node (J) at (0,6) {$\mathbf{(3,0)}$};

    \path [->] (A) edge (B);
    \path [->] (B) edge (C);
    \path [->] (B) edge (D);
    \path [->] (C) edge (E);
    \path [->] (D) edge (E);
    \path [->] (D) edge (F);
    \path [->] (E) edge (G);
    \path [->] (E) edge (H);
    \path [->] (F) edge (H);
    \path [->] (G) edge (I);
    \path [->] (H) edge (I);
    \path [->] (I) edge (J);
\end{tikzpicture}
\par}
\vspace{1cm}
{\centering
\begin{tikzpicture}[scale=0.8]
     {\color{black!50!white}
    \node (C-5) at (3,-5) {$(-5,0)$};
    \node (C-4) at (3,-4) {$(-4,1)$};
    \node (B-3) at (2,-3) {$(-3,0)$};
    \node (D-3) at (4,-3) {$(-3,2)$};
    \node (B-2) at (2,-2) {$(-2,1)$};
    \node (D-2) at (4,-2) {$(-2,3)$};
    \node (A-1) at (1,-1) {$(-1,0)$};
    \node (C-1) at (3,-1) {$(-1,2)$};
    \node (E-1) at (5,-1) {$(-1,4)$};
    \node (A0) at (1,0) {$(0,1)$};
    \node (C0) at (3,0) {$(0,3)$};
    \node (E0) at (5,0) {$(0,5)$}; }
    \node (A1) at (1,1) {$\mathbf{(1,0)}$};
    \node (C1) at (3,1) { $(1,2)$ };
    \node (E1) at (5,1) { $(1,4)$ };
    \node (B2) at (2,2) { $\mathbf{(2,1)}$ };
    \node (D2) at (4,2) { $(2,3)$ };
    \node (B3) at (2,3) { $\mathbf{(3,0)}$ };
    \node (D3) at (4,3) { $\mathbf{(3,2)}$ };
    \node (C4) at (3,4) { $\mathbf{(4,1)}$ };
    \node (C5) at (3,5) { $\mathbf{(5,0)}$ }; 
    
    \path [->] (C-5) edge (C-4);

    \path [->] (C-4) edge (B-3);
    \path [->] (C-4) edge (D-3);
    
    \path [->] (B-3) edge (B-2);
    \path [->] (D-3) edge (B-2);
    \path [->] (D-3) edge (D-2);
    
    \path [->] (B-2) edge (A-1);
    \path [->] (B-2) edge (C-1);
    \path [->] (D-2) edge (C-1);
    \path [->] (D-2) edge (E-1);
    
    \path [->] (A-1) edge (A0);
    \path [->] (C-1) edge (A0);
    \path [->] (C-1) edge (C0);
    \path [->] (E-1) edge (C0);
    \path [->] (E-1) edge (E0);
    
    \path [->] (A0) edge (A1);
    \path [->] (A0) edge (C1);
    \path [->] (C0) edge (C1);
    \path [->] (C0) edge (E1);
    \path [->] (E0) edge (E1);

    \path [->] (A1) edge (B2);
    \path [->] (C1) edge (B2);
    \path [->] (C1) edge (D2);
    \path [->] (E1) edge (D2);

    \path [->] (B2) edge (B3);
    \path [->] (B2) edge (D3);
    \path [->] (D2) edge (D3);

    \path [->] (B3) edge (C4);
    \path [->] (D3) edge (C4);

    \path [->] (C4) edge (C5);       
\end{tikzpicture}
\par}
\end{minipage}
\begin{minipage}[b][0.6\textheight][c]
{0.5\textwidth}  
\begin{tikzpicture}[scale=0.8]
    {\color{black!50!white}
    \node (-7 0) at (4,-7) {$(-7,0)$};

    \node (-6 1) at (4,-6) {$(-6,1)$};

    \node (-5 0) at (3,-5) {$(-5,0)$};
    \node (-5 2) at (5,-5) {$(-5,2)$};

    \node (-4 1) at (3,-4) {$(-4,1)$};
    \node (-4 3) at (5,-4) {$(-4,3)$};

    \node (-3 0) at (2,-3) {$(-3,0)$};
    \node (-3 2) at (4,-3) {$(-3,2)$};
    \node (-3 4) at (6,-3) {$(-3,4)$};

    \node (-2 1) at (2,-2) {$(-2,1)$};
    \node (-2 3) at (4,-2) {$(-2,3)$};
    \node (-2 5) at (6,-2) {$(-2,5)$};
    
    \node (-1 0) at (1,-1) {$(-1,0)$};
    \node (-1 2) at (3,-1) {$(-1,2)$};
    \node (-1 4) at (5,-1) {$(-1,4)$};
    \node (-1 6) at (7,-1) {$(-1,6)$};

    \node (0 1) at (1,0) {$(0,1)$};
    \node (0 3) at (3,0) {$(0,3)$};
    \node (0 5) at (5,0) {$(0,5)$};
    \node (0 7) at (7,0) {$(0,7)$};
    }
    
    \node (1 0) at (1,1) {$\mathbf{(1,0)}$};
    \node (1 2) at (3,1) {$(1,2)$};
    \node (1 4) at (5,1) {$(1,4)$};
    \node (1 6) at (7,1) {$(1,6)$};
    
    \node (2 1) at (2,2) {$\mathbf{(2,1)}$};
    \node (2 3) at (4,2) {$(2,3)$};
    \node (2 5) at (6,2) {$(2,5)$};
    
    \node (3 0) at (2,3) {$\mathbf{(3,0)}$};
    \node (3 2) at (4,3) {$\mathbf{(3,2)}$};
    \node (3 4) at (6,3) {$(3,4)$};
    
    \node (4 1) at (3,4) {$\mathbf{(4,1)}$};
    \node (4 3) at (5,4) {$\mathbf{(4,3)}$};

    \node (5 0) at (3,5) {$\mathbf{(5,0)}$};
    \node (5 2) at (5,5) {$\mathbf{(5,2)}$};

    \node (6 1) at (4,6) {$\mathbf{(6,1)}$};

    \node (7 0) at (4,7) {$\mathbf{(7,0)}$};

    \path [->] (-7 0) edge (-6 1);

    \path [->] (-6 1) edge (-5 0);
    \path [->] (-6 1) edge (-5 2);

    \path [->] (-5 0) edge (-4 1);
    \path [->] (-5 2) edge (-4 1);
    \path [->] (-5 2) edge (-4 3);

    \path [->] (-4 1) edge (-3 0);
    \path [->] (-4 1) edge (-3 2);
    \path [->] (-4 3) edge (-3 2);
    \path [->] (-4 3) edge (-3 4);

    \path [->] (-3 0) edge (-2 1);
    \path [->] (-3 2) edge (-2 1);
    \path [->] (-3 2) edge (-2 3);
    \path [->] (-3 4) edge (-2 3);
    \path [->] (-3 4) edge (-2 5);

    \path [->] (-2 1) edge (-1 0);
    \path [->] (-2 1) edge (-1 2);
    \path [->] (-2 3) edge (-1 2);
    \path [->] (-2 3) edge (-1 4);
    \path [->] (-2 5) edge (-1 4);
    \path [->] (-2 5) edge (-1 6);

    \path [->] (-1 0) edge (0 1);
    \path [->] (-1 2) edge (0 1);
    \path [->] (-1 2) edge (0 3);
    \path [->] (-1 4) edge (0 3);
    \path [->] (-1 4) edge (0 5);
    \path [->] (-1 6) edge (0 5);
    \path [->] (-1 6) edge (0 7);

    \path [->] (0 1) edge (1 0);
    \path [->] (0 1) edge (1 2);
    \path [->] (0 3) edge (1 2);
    \path [->] (0 3) edge (1 4);
    \path [->] (0 5) edge (1 4);
    \path [->] (0 5) edge (1 6);
    \path [->] (0 7) edge (1 6);

    \path [->] (1 0) edge (2 1);
    \path [->] (1 2) edge (2 1);
    \path [->] (1 2) edge (2 3);
    \path [->] (1 4) edge (2 3);
    \path [->] (1 4) edge (2 5);
    \path [->] (1 6) edge (2 5);

    \path [->] (2 1) edge (3 0);
    \path [->] (2 1) edge (3 2);
    \path [->] (2 3) edge (3 2);
    \path [->] (2 3) edge (3 4);
    \path [->] (2 5) edge (3 4);

    \path [->] (3 0) edge (4 1);
    \path [->] (3 2) edge (4 1);
    \path [->] (3 2) edge (4 3);
    \path [->] (3 4) edge (4 3);

    \path [->] (4 1) edge (5 0);
    \path [->] (4 1) edge (5 2);
    \path [->] (4 3) edge (5 2);

    \path [->] (5 0) edge (6 1);
    \path [->] (5 2) edge (6 1);

    \path [->] (6 1) edge (7 0);    
\end{tikzpicture}
\end{minipage}
\par}
\caption{Hasse diagrams for $n=3,5,7$, in the quotient poset of $(\rho,\alpha)$-tables.
In boldface, tables of type \textbf{b}, in normal type those of type \textbf{c},
and greyed out those of type 
\textbf{a}.\label{fig:ratables}}
\end{figure}
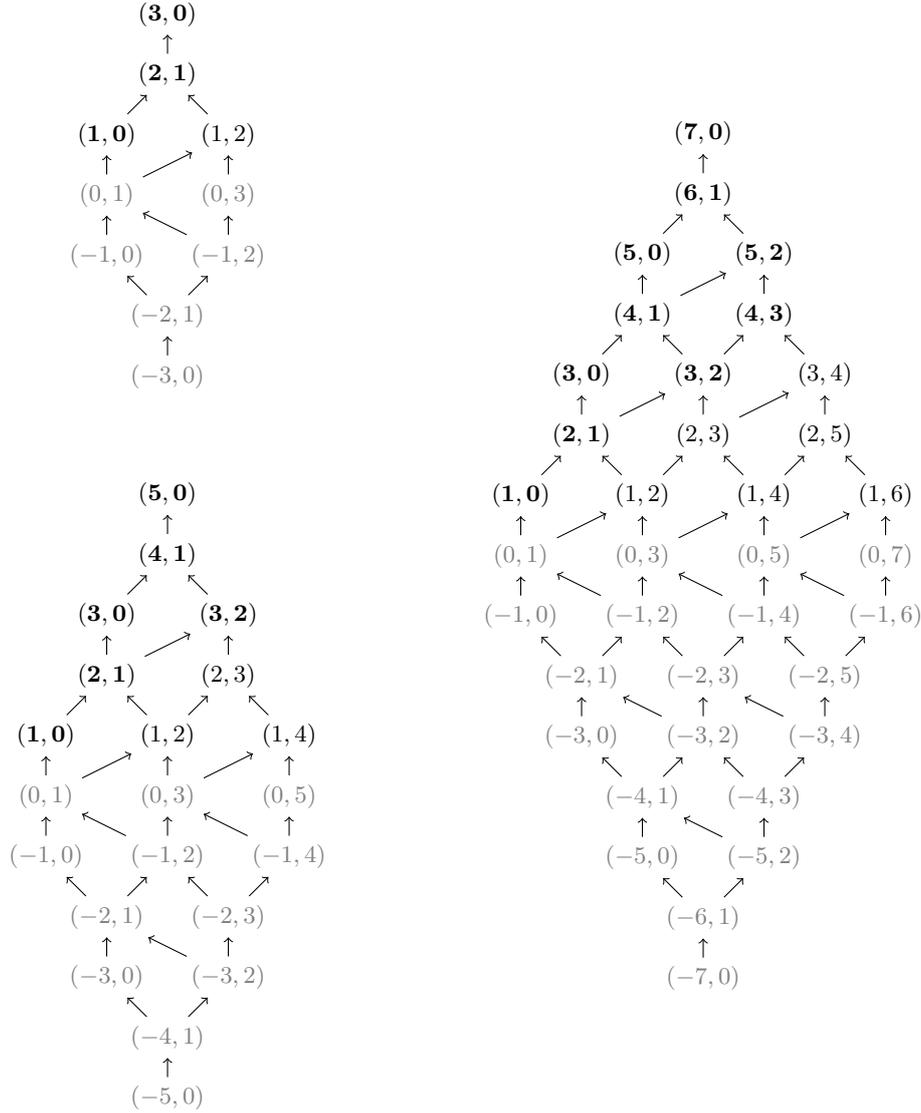

\begin{example}
Figure \ref{fig:ratables} illustrates the resulting poset of $(\rho, \alpha)$-tables,
for $n=3,5,7$, with different intensities according to types 
\textbf{a}, \textbf{b}, \textbf{c} of Theorem \ref{thm:optrule}.
\qed
\end{example}

The relations (\ref{eq:transreduc}) defining the transitive reduction in $(\mathbb T,\le)$ translate to
\begin{equation}\label{eq:transreduc2}
\begin{split}
(\rho, \alpha)&\le(\rho+1, \alpha-1) \\
(\rho, \alpha)&\le(\rho+1, \alpha+1) \ 
\end{split} 
\end{equation}
in the quotient poset, and it is easy to see, while proving that Lemma \ref{lem:preorder2poset}
is applicable to $(\mathbb T,\le)$ (see the appendix),
that (\ref{eq:transreduc2}) is precisely the transitive reduction in the quotient poset. 
This makes the Hasse diagrams 
like those of Figure \ref{fig:ratables} very easy to generate for any $n$. 

\bigskip
The classical premiss-based rule coincides with the one formed exactly by the tables of type \textbf{b}
(see Example \ref{ex:classicalRules}). This suggests that it is possible to characterise exactly under what
conditions on $w$ and $\theta$ the premiss-based rule is the optimal one. The result is 
given in Theorem \ref{thm:pb-optimal}. It will be a consequence of the following lemma.
In the sequel we will denote $G_{\rho,\alpha}$ the function defined in (\ref{eq:G}).
\begin{lemma}\label{lem:pb-optimal}
In the poset $(\mathbb T/\negthickspace\sim,\le)$ of $(\rho,\alpha)$-tables, 

\begin{enumerate}
\item The subset of tables of type \textbf{b} has a unique minimal element: Table $(1,0)$.
\item The union of tables of type \textbf{a} and \textbf{c} has a unique maximal element: 
Table $(\frac{n-1}{2},\frac{n+1}{2})$. 
\end{enumerate}
\end{lemma} 
\begin{proof}
The statements are equivalent to say that $G_{\rho,\alpha}\le G_{1,0}$ for $\rho>\alpha\ge 0$, 
and that $G_{\rho, \alpha}\ge G_{\frac{n-1}{2},\frac{n+1}{2}}$ 
for $\rho<\alpha$ with $\alpha\ge 0$. For the first inequality, both exponents $-\rho-\alpha$ and $-\rho+\alpha$
are negative, hence $G_{\rho, \alpha}(\eta)\le\eta^{-1}+\eta^{-1} = G_{1,0}(\eta)$; for the second,
the first exponent is positive, and the second is greater or equal to $-n$, hence 
$G_{\rho, \alpha}(\eta)\ge\eta+\eta^{-n} = G_{\frac{n-1}{2},\frac{n+1}{2}}(\eta)$.
\end{proof}

\begin{theorem}\label{thm:pb-optimal}
Let $0<w<1$ and $\frac{1}{2}<\theta<1$ be the given weight and competence level, and 
$n$ the committee size. The premiss-based rule is optimal if and only if 
\begin{equation}
\theta\ge w
\quad\text{ and }\quad
\frac{\theta}{1-\theta}+\Big(\frac{\theta}{1-\theta}\Big)^{-n}\ge \frac{2(1-w)}{w}
\label{eq:pb-optimal}
\end{equation}
\end{theorem}
\begin{proof}
Denote, as before, $\eta:=\frac{\theta}{1-\theta}$, and set also $\xi:=\frac{2(1-w)}{w}$.
In view of Lemma \ref{lem:pb-optimal}, the necessary and sufficient condition for the
premiss-based rule to be optimal is that the point $(\eta,\xi)$ lie above the curve 
$G_{1,0}$ and below the curve $G_{\frac{n-1}{2},\frac{n+1}{2}}$. That means 
$2\eta^{-1}\le\xi\le \eta+\eta^{-n}$, or equivalently,
\begin{equation*}
2\Big(\frac{\theta}{1-\theta}\Big)^{-1}\le
\frac{2(1-w)}{w}\le 
\frac{\theta}{1-\theta}+\Big(\frac{\theta}{1-\theta}\Big)^{-n} 
\ .
\end{equation*}
The first inequality is equivalent to $\theta\ge w$, and we are done.
\end{proof}
 
\begin{figure}[h]
{\centering\small
\usetikzlibrary{calc}
\begin{tikzpicture}[yscale=5, xscale=10]  
\begin{scope}[shift={(-0.4,0)}]
  \draw[->] (0.5,0) -- (1.1,0);
  \draw[-, dashed] (0.4,0) -- (0.5,0);
  \node[below, align=center] at (1.1,-0.02) {$\theta$};
  \draw[->] (0.4,0) -- (0.4,1.2);
  \node[left, align=center] at (0.39,1.2) {$w$};
  \draw[-] (.39,1) -- (0.41,1); 
  \node[left] at (0.39,1) {$1$};
  \draw[-] (.39,0) -- (0.41,0); 
  \node[left] at (0.39,0) {$0$};
  \draw[-] (0.39,0.5858) -- (0.41,0.5858); 
  \node[left] at (0.39,0.5858) {\tiny$2-\sqrt{2}$};
  \draw[-] (0.5,-0.02) -- (0.5,0.02); 
  \node[below, align=center] at (0.5,-0.02) {$\frac{1}{2}$};
  \draw[-] (1,-0.02) -- (1,0.02); 
  \node[below, align=center] at (1,-0.02) {$1$};
  \draw[-, dashed] (0.4,0.5858) -- (0.5858,0.5858);
  \draw[-, dashed] (0.5858, 0) -- (0.5858,0.5858);
  \draw[-] (0.5858,-0.02) -- (0.5858,0.02); 
  \node[below, align=center] at (0.5858,-0.02) {\tiny$2-\sqrt{2}$};
    \draw[thick, x=1cm, y=1cm] 
    plot[domain=0.5:1, samples=144, smooth] (\x, \x);
    \draw[thick, x=1cm, y=1cm] 
    plot[domain=0.5:1, samples=144, smooth] (\x, {(2-2*\x)/(2-\x)});
    \draw[x=1cm, y=1cm] 
    plot[domain=0.5:1, samples=144, smooth] (\x, {(2*(1-\x))/(\x+((1-\x)^(3+1))/(\x^3)+2*(1-\x))});
    \draw[x=1cm, y=1cm] 
    plot[domain=0.5:1, samples=144, smooth] (\x, {(2*(1-\x))/(\x+((1-\x)^(5+1))/(\x^5)+2*(1-\x))});
    \draw[x=1cm, y=1cm] 
    plot[domain=0.5:1, samples=144, smooth] (\x, {(2*(1-\x))/(\x+((1-\x)^(7+1))/(\x^7)+2*(1-\x))});
  \node[above right] at (0.92,0.1) {$\theta=\frac{2-2w}{2-w}$};
  \node[above right] at (0.9,0.8) {$\theta=w$};
  \node[right] at (0.77,0.55) {\footnotesize pb-optimality};  
  \draw[-, color=black!50!white, dotted] (0.5,1) -- (1,1);
  \draw[-, color=black!50!white, dotted] (1,0) -- (1,1);
  \draw[-, color=black!50!white, dotted] (0.5,0) -- (0.5,1);  
\end{scope}
\end{tikzpicture}
\caption{Region of optimality of the premiss-based (pb) rule, in the natural 
coordinates $(\theta,w)$. The thinner curves correspond to $n=3,5,7$ approaching 
monotonically the curve $\theta=\frac{2-2w}{2-w}$ as $n\to\infty$. 
If $\theta<w$, some tables of type \textbf{b} must leave the set ${\mathbb T}^+$; 
if $(\theta,w)$ falls below the lower curve, other tables must join those of 
type \textbf{b} in the set ${\mathbb T}^+$; and both things may happen at the same time, 
for $\theta<2-\sqrt{2}$.}
\label{fig:pb-opt}
}
\end{figure}
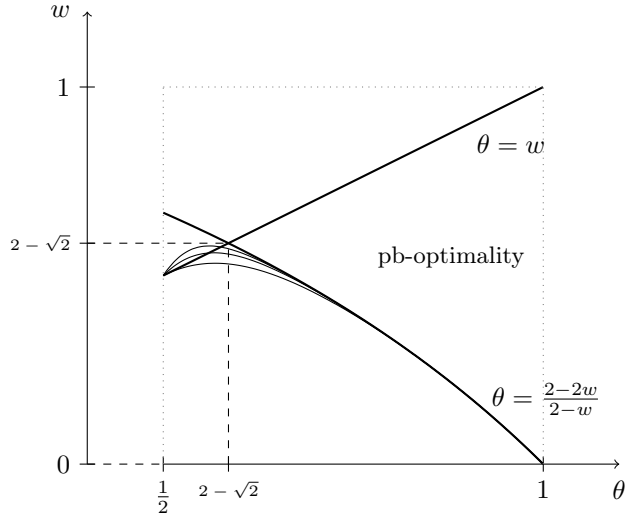

A simpler sufficient condition, independent of the committee size, is given in the next corollary.
Figure \ref{fig:pb-opt} illustrates both theorem and corollary.
\begin{corollary}\label{cor:pb-optimal}
Let $0<w<1$ and $\frac{1}{2}<\theta<1$ be as in Theorem \ref{thm:pb-optimal}. If
\begin{equation*}
\theta\ge w 
\quad\text{ and }\quad
\theta \ge \frac{2(1-w)}{2-w}
\ ,
\end{equation*}
then the premiss-based rule is optimal, for all committee sizes.
\end{corollary} 
\begin{proof}
The second inequality results from ignoring the negative exponential term in (\ref{eq:pb-optimal}). 
\end{proof}

\begin{example}
From the corollary, in the balanced case $w=\frac{1}{2}$, 
one can be sure that the premiss-based rule is optimal as soon as the
competence level is above $\frac{2}{3}$. 
If this is not the case, other voting tables (those of type \textbf{c}) 
have to be successively added to the set of positive tables, as
the competence level decreases. Table \ref{tab:theta0} shows the 
critical $\theta_0$, to four decimal places, for the tables of type \textbf{c} with $n\le 7$.
Notice that the value of $\theta_0$ is independent of $n$. 
\begin{table}[h] \centering
    {\extrarowheight 2pt \footnotesize
  \begin{tabular}{c|c|c|c|c|c|c}
     $(\rho,\alpha)$ & $(3,4)$ & $(2,3)$ & $(1,2)$ & $(2,5)$ & $(1,4)$ & $(1,6)$  \\
    \hline 
    $\theta_0$ & $0.6658$ & $0.6628$ & $0.6478$ & $0.5449$ & $0.5326$ & $0.5141$ 
  \end{tabular}
}
\caption{Positive tables of type \textbf{c} when $w=\frac{1}{2}$ and 
$\theta\in(\frac{1}{2},\theta_0)$.}\label{tab:theta0}
\end{table}  

One might conjecture that the tables of type \textbf{c} enter the optimal rule
following the increasing value of $-\rho+\alpha$, an among those with the same value, 
following the decreasing value of $\rho+\alpha$. This is true up to $n=11$. For $n=13$ this 
regularity breaks down and table $(3,10)$ enters before $(1,6)$, at $\theta_0=0.5160$.
Hence, we do not find here any computational shortcut to determine the optimal rule
for general $n$, even in the case $w=\frac{1}{2}$.
\qed  
\end{example}

By contrast, the classical conclusion-based rule is \emph{never} optimal: For any $n$, a pair of tables $(x,y,z,t)$ 
and $(x',y',z',t')$ can be found that lead to different results according to conclusion-based, and yet they
belong to the same $(\rho,\alpha)$ class. For instance, for $n=3$, we have $(2,0,0,1)$ leading to $C$ and 
$(1,1,1,0)$ leading to $\neg C$; but both belong to the class $(1,0)$ and should have the same status
in the optimal rule. 
However, it is worth noting that the conclusion-based rule can be better than the premiss-based 
for certain values of $w$ and $\theta$.

\section{Computations and software}\label{sec:ComputSoft}
The optimality condition adds
some more relations to the partial order defined by the admissibility requirement. 
This simplifies further the computation of the optimal rule.

\begin{proposition}\label{pr:OptReduction}
If $(\rho,\alpha)$ is a positive table in the optimal rule, then $(\rho,\alpha-2)$ is 
also a positive table in the optimal rule, whenever these values make sense.

Adding the relations $(\rho,\alpha)\le(\rho,\alpha-2)$ to $(\mathbb T,\le)$, the new
partial order has the transitive reduction defined by
\begin{equation*}
\begin{split}
(\rho, \alpha)&\le(\rho, \alpha-2) \\
(\rho, \alpha)&\le(\rho+1, \alpha+1) \\
(n-1,1)&\le(n,0)
\end{split} 
\end{equation*}
\end{proposition}

\begin{proof}
We need to see that for the functions $G_{\rho,\alpha}$ defined in (\ref{eq:G}),
we have $G_{\rho,\alpha-2}\le G_{\rho,\alpha}$ in their whole domain $(1,\infty)$. 
The inequality can be expressed as
\begin{equation*}
0\le \eta^{2\alpha}-\eta^{2\alpha-2}-\eta^2+1 = (\eta^{2\alpha-2}-1)(\eta^2-1)
\ ,
\end{equation*}
which is obviously true for all $\eta>1$ and $\alpha\ge 2$. 

The first relation in the transitive reduction (\ref{eq:transreduc2}) 
is no longer present in the new transitive reduction, since now there is
an element in between:
\begin{equation*}
(\rho, \alpha)\le(\rho, \alpha-2)\le(\rho+1, \alpha-1) 
\end{equation*}
except in $(n-1,1)\le(n,0)$. 
\end{proof}

\begin{figure}
{\centering\footnotesize
\begin{tikzpicture}
    \node (A) at (0,0) {$(-3,0)$};
    \node (B) at (0,1) {$(-2,1)$};
    \node (C) at (0,2) {$(-1,2)$};
    \node (D) at (-1,3) {$(-1,0)$};
    \node (E) at (1,3) {$(0,3)$};
    \node (F) at (0,4) {$(0,1)$} ;
    \node (G) at (0,5) {$(1,2)$} ;
    \node (H) at (0,6) {$(1,0)$} ;
    \node (I) at (0,7) {$(2,1)$} ;
    \node (J) at (0,8) {$(3,0)$} ;

    \path [->] (A) edge (B);
    \path [->] (B) edge (C);
    \path [->] (C) edge (D);
    \path [->] (C) edge (E);
    \path [->] (D) edge (F);
    \path [->] (E) edge (F);
    \path [->] (F) edge (G);
    \path [->] (G) edge (H);
    \path [->] (H) edge (I);
    \path [->] (I) edge (J);
\end{tikzpicture}
\par
}
\caption{Partial order in $\mathbb T$ induced by the optimality condition (see Proposition \ref{pr:OptReduction}). \label{fig:Hasse-opt}}
\end{figure}
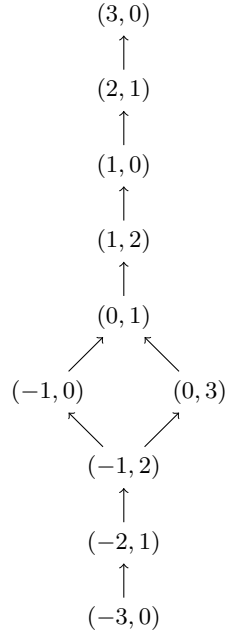

The resulting Hasse diagram is ``thinner'',
and the total number of upper sets is reduced.
As an example, the case $n=3$ is depicted in Figure \ref{fig:Hasse-opt}. 
There are only twelve upper sets left 
after this simplification. The new poset is still ranked, but 
$\rho=x-t$ is no longer a rank function.

This reduction is relevant if we are only interested in the optimal rule.
If we want to build instead a ranking of rules, then 
Proposition \ref{pr:OptReduction} is not useful.

We have built a program in Python, with a graphical interface, that allows the user
to input the values of $n$, $w$ and $\theta$, and produces a ranking of decision rules.
It can be currently found as a public Mercurial repository in \url{https://discursive-dilemma.sourceforge.io/}, 
or requested 
directly to the authors. The program allows to specify different competence levels for
the different committee members (an extension discussed in the next section), so that
in fact formulae (\ref{eq:PFP}) and (\ref{eq:PFN}) are not used, but instead 
the probability to get a voting table $(x,y,z,t)$ is computed taking into account
all possible permutations of voters.

\begin{figure}
\centering
\includegraphics{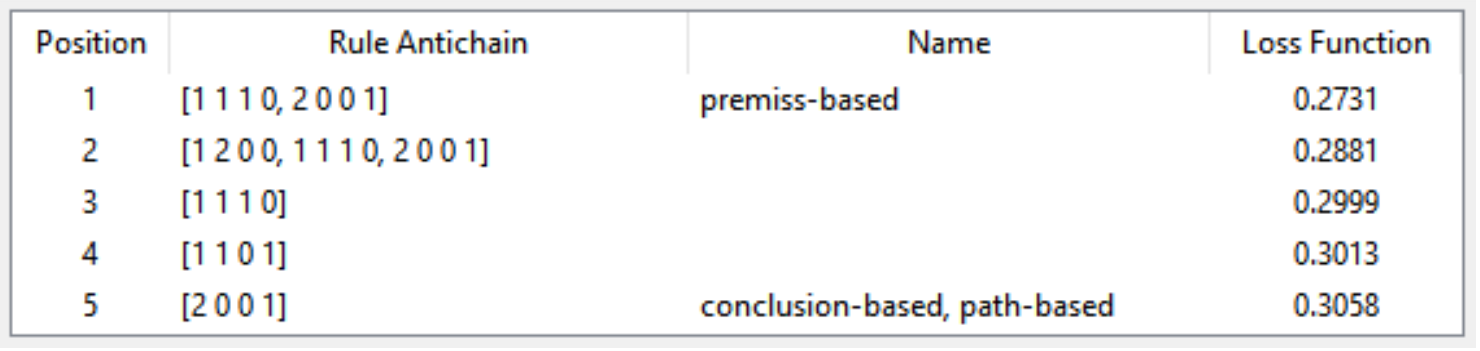}\\
\vspace{2mm}
\includegraphics{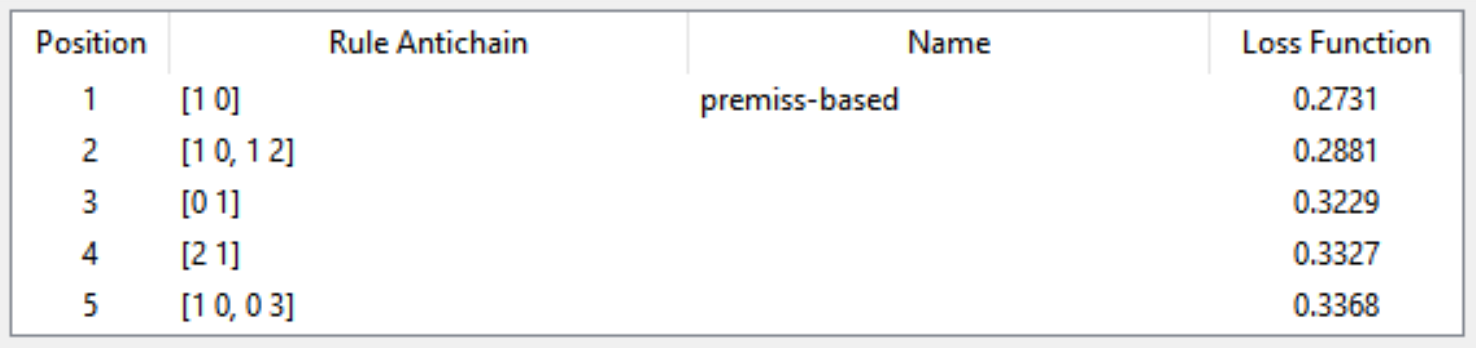}
\caption{The two rankings produced by the Python code: For rules in the form $(x,y,z,t)$ 
and for rules in the form $(\rho, \alpha)$. In this example, the input was 
$n=3$, $w=0.5$, $\theta=(0.6,0.7,0.8)$.}
\label{fig:output}
\end{figure}

As an example of output, see Figure \ref{fig:output}. Two rankings are produced, the first 
corresponding to voting tables and rules in extended form $(x,y,z,t)$, and the second 
in compact form $(\rho, \alpha)$. They are not a direct translation of each other, since
a rule that cannot be expressed in compact form (because members of the same $(\rho,\alpha)$ 
class are assigned different conclusions) may actually be better than the next rule 
respecting the equivalence relation. 

The rules are expressed by means of the antichain that determines the upper set of 
positive tables. A name is printed if the rule is one of premiss-based, conclusion-based
or path-based, and the value of the loss function (\ref{eq:WAOT}) of each rule is also given. 
Notice that in this example the rules
in positions 3 to 5 in the extended version are not expressible in compact form, but
are better than the third rule in the second ranking. Of course, the optimal rule will
always coincide in both rankings. 

At the moment, the program only outputs up to the five best rules, but this is an arbitrary 
parameter that can be easily changed in the source code. Also, we have not made
any special effort for efficiency. It has been conceived only as a playground and 
checking tool. 

\bigskip

Concerning the computational complexity of producing the admissible decision rules, let us note
first that the total number of voting tables is equal to 
$\frac{1}{24}(2n^3+15n^2+34n+21)$ 
in the original poset $(\mathbb T,\le)$, after identifying transposed tables. 
In the quotient poset  $(\mathbb T/\negthickspace\sim,\le)$, this number is reduced to 
$\frac{1}{2}(n+1)(n+2)$. 

Since admissible rules are in bijection with upper sets, and these in turn 
are determined by their antichains, 
we first identify the latter. For this 
task we make use of the Python package \verb|networkx|, which contains a function \verb|antichains()|.
The generation of the corresponding upper set and the evaluation of each table 
contained in it is very easy. 
The contribution to the loss function of tables previously computed is stored
to speed up the computations. 
The maximum cardinality of an antichain (after identifying transposed tables) 
is $\frac{1}{8}(n+3)(n+1)$ in $(\mathbb T,\le)$ and 
$\frac{1}{2}(n+1)$ in $(\mathbb T/\negthickspace\sim,\le)$.
We sketch in the appendix the computation of these numbers.
 
\bigskip
Apart from having the optimal rule or a ranking of the best rules, one might be just interested
in knowing which conclusion has to be assigned to a given voting table $T$ under the optimal rule. 
This is very easy by asserting inequality (\ref{eq:goodTable}), in the case of equal competences.
Otherwise, the probabilities of a False Positive and a False Negative have
to be computed for that table, taking into account the different competences, and then determine  
their contribution to the loss function.

Finally, in the equal competences case, one may like to determine, given a fixed
weight $w$, the intervals of
competence $\theta$ where $r(T)=0$ or $1$ for the optimal rule $r$. We 
need to find the root or roots of the equation 
\begin{equation}\label{eq:findroot}
G_{\rho,\alpha}(\eta)=\eta^{-\rho-\alpha}+\eta^{-\rho+\alpha}=\frac{2(1-w)}{w}
\ .
\end{equation}
This is very easy numerically. The functions involved are simple
to evaluate so that pure bisection, for instance, is very fast. We only need to have a bracket 
where the roots are guaranteed to lie. Indeed, they are readily found:
\begin{description}
\item [Tables of type a:] ($\rho\le 0$). \\
If $w\ge \frac{1}{2}$ the table is bad.\\
If $w<\frac{1}{2}$, the unique root of (\ref{eq:findroot}) is less than the
root of $\eta^{-\rho+\alpha}=\frac{2(1-w)}{w}$. Therefore, the solution
to (\ref{eq:findroot}) 
will be found between 1 and
\begin{equation}\label{eq:bound_a}
\Big(\frac{2(1-w)}{w}\Big)^{\frac{1}{-\rho+\alpha}}\ .
\end{equation}
\item [Tables of type b:] ($\rho\ge\alpha$). \\
If $w\le \frac{1}{2}$ the table is good.\\
If $w>\frac{1}{2}$, the unique root of (\ref{eq:findroot}) is less than the
root of $2\eta^{-\rho+\alpha}=\frac{2(1-w)}{w}$, and it will be found
between 1 and 
\begin{equation*}
\Big(\frac{w}{1-w}\Big)^{\frac{1}{\rho-\alpha}}\ .
\end{equation*}
\item [Tables of type c:] ($0<\rho<\alpha$). \\
If $w\le \frac{1}{2}$ there is a unique root, and the same bound (\ref{eq:bound_a})
is valid.\\
If $\frac{1}{2}<w\le\frac{2}{3}$, we first check if the minimum $\eta^*$ of $G_{\rho,\alpha}$ (see (\ref{eq:minG}))
satisfies $G_{\rho,\alpha}(\eta^*)>\frac{2(1-w)}{w}$, in which case there are no roots of (\ref{eq:findroot}) and
the table is bad; or $G_{\rho,\alpha}(\eta^*)<\frac{2(1-w)}{w}$ and there are two roots: The first between
1 and $\eta^*$, and the second between $\eta^*$ and (\ref{eq:bound_a}).\\
If $w>\frac{2}{3}$, the table is bad, because then $\frac{2(1-w)}{w}<1$, whereas $G_{\rho,\alpha}(\eta^*)>1$ always.
\end{description}

\section{Conclusions and discussion}\label{sec:ConcDisc}
We have studied the discursive dilemma in its simplest classical form, and proposed 
a method to obtain the best rule (or a ranking of rules) by minimising a loss function
that combines false positives and false negatives.
Actually, we have introduced a
family of loss functions, parametrised by the number $0<w<1$.

The decision rules considered satisfy very mild and reasonable conditions of symmetry 
and monotonicity (Definition \ref{def:admissible}). In fact,
the second condition is not necessary a priori if one is only
interested in the best rule and not in ranking rules. In that case, monotonicity
appears a posteriori as a property of the optimal rule. 

Generically, the optimal rule will be unique, but specific values of weight $w$
and competence $\theta$ may lead to ties in the evaluation of the loss function,
in particular in its minimum value. To make the exposition simpler, we have avoided
mentioning this possibility throughout the paper. 

The loss function is a modelling choice. In any real instance it must be chosen
to reflect what the best rule is intended to achieve. 
The important point is that the optimisation setting is 
worth considering for problems of judgment aggregation in general. 

In Alabert-Farre \cite{AlabertFarre2022}, where this point of view was introduced, 
some possible extensions and open problems were
discussed substantially. We summarise them here:

\begin{itemize}
\item \emph{Different competence for each voter.} 
This is the simplest extension. If $J_k$ is the 
voting table consisting only on the vote of voter $k$, with competence level $\theta_k$, 
the resulting table is $J_1+\dots+J_n$, whose probability law can also be computed,
and the probabilities of false positive and false negative will be 
\begin{align*}
\mathbb P_r(\text{FP})
  &= 
  \sum_{\{r(x,y,z,t)=1\}} \mathbb P_{P\wedge\neg Q}\{J_1+\cdots+J_n=(x,y,z,t)\} \ ,  
  \\  
\mathbb P_r(\text{FN}) 
  &= \sum_{\{r(x,y,z,t)=0\}} \mathbb P_{P\wedge Q}\{J_1+\cdots+J_n=(x,y,z,t)\} 
  \ .  
\end{align*}
Our software already computes the ranking of rules in this more general situation,
\item \emph{Different competence for each premiss or state of nature.} 
The competence of a voter may be in fact a vector 
$\theta=(\theta_P, \theta_{\neg P}, \theta_Q, \theta_{\neg Q})$ of competences
depending on the premiss and/or the true state of nature. In List \cite{List2005}, 
the probability of appearance of the doctrinal paradox is studied also when
the competence is different on $P$ and $Q$. 
The computation of $\mathbb P_r(\text{FP})$ and $\mathbb P_r(\text{FN})$ is more
involved in this case, but still feasible.

\item \emph{Non-independence between voters.} 
If the committee members do not vote independently, perhaps through a deliberation
process with influential individuals, then the full joint law of the vector
$(J_1,\dots, J_n)$ of individual voting tables is needed to compute the law
of the sum $J_1+\cdots+J_n$ under the different states of nature. 
Boland \cite{10.2307/2348873} studied this situation for the voting of a single
question assuming the presence of a ``leader'' in the committee. Other works
that studied epistemic social choice with correlated voters in the last decade
include Peleg and Zamir \cite{Peleg2012}, Dietrich and Spiekermann \cite{dietrich_spiekermann_2013}
\cite{10.1093/mind/fzt074}, and Pivato \cite{PIVATO201751}.

\item \emph{Non-independence between premisses.} In practical examples, the 
premisses $P$ and $Q$ can very well be interconnected, in the sense that believing that
$P$ is true or false can change the perception on the truth or falsity of $Q$.
This may lead to a different competence in asserting $Q$ depending on the
decision on $P$. Then the joint law of the competences under the four states of 
nature are needed to complete the computations. 

The extreme case where one combination of premisses is impossible it is 
treated in Bozbay \cite{Bozbay2019}, where in addition abstentions are allowed.

\item \emph{More than two premisses.} There is no difficulty in extending the
setting to a conjunctive agenda with any number of premisses $P_1,\dots, P_s$.
A voting table will be an element of 
$\mathbb T=\{(x_1,\dots,x_{2^s})\in \mathbb N^{2^s}:\ \sum_{i=1}^{2^s} x_i=n\}$. 
The concepts of admissible rule and of false positive are easily extended.
\end{itemize}

Note that \emph{disjunctive agendas}, in which the conclusion is true if and only if at least one premiss
is true, are dual to the conjunctive case, by negation of the doctrine (see 
List \cite{List2005}, Bovens and Rabinowicz \cite{Bovens2006}, or Miyashita \cite{Miyashita2021}).
They can be considered easily within our framework.
  
An extension with an obvious practical interest is allowing abstentions, or committees
with an even number of members. It is clear that enforcing an opinion on all clauses
of the agenda may be inconvenient or simply impossible. These so-called \emph{incomplete judgments}
have been considered in Gärdenfors \cite{Gardenfors2006}, Dietrich-List \cite{DietrichList2008}, Terzopoulou and Endriss \cite{TerzopoulouEndrissSAGT2019},
and Bozbay \cite{Bozbay2019}.

It is natural to ask which desirable properties satisfies the optimal rule of a given criterion. 
We leave this as an open question. In relation to the classical axioms of judgment
aggregation and their (im)possibility theorems (see e.g.~List \cite{List2012}), 
and since here we are centred in reaching
a right conclusion for whatever reasons, collective rationality can only be achieved by assigning a value
to the premisses \emph{after} deciding on the conclusion (see Pigozzi et al. \cite{Pigozzi2009});
but then the properties of monotonicity (in the classical sense), unanimity and systematicity 
need not be satisfied on the whole agenda. 
On the other hand, the anonymity 
requirement is trivially met in our setting.        
In any case, the advantage of the optimisation model is the immediate existence 
of decision rules; each of the rules is evaluated through a real-valued loss function, 
hence at least one rule with a minimal value must exist. 
Distance-based methods to reach consensus share this feature.
    
\section*{Appendix}\label{sec:Appendix} 
$\bullet$ \emph{Proof of the equivalence between admissibility and upper sets (Proposition \ref{pr:goodUpper}):}

Given $r \in \mathcal{A}$, the set \{$u \in \mathbb{T} : r(u) = 1$\} is an upper set of $\mathbb{T}$.
We want to see that for all $x \in S = \{u \in \mathbb{T} : r(u) = 1\}$ and for all $y \in \mathbb{T}$, 
such that $x < y$, we must have $y \in S$. Take $x \in S$ and $y \in \mathbb{T}$ such that $x < y$. 
Since $r \in \mathcal{A}$ we have $1 = r(x) < r(y)$, hence $r(y) = 1$, and therefore $y \in S$.

$\bullet$ \emph{Proof of the final claim in Subsection \ref{ssec:OptCrit}:}

In case $\theta=\frac{1}{2}$, the problem is trivial because denoting 
$\Sigma=\sum\limits_{r(x,y,z,t)=0} \frac{n!}{x!y!z!t!} (\frac{1}{2})^{2(x+y+z+t)}$, 
we have
$L_w(r)= w(1-\Sigma)+(1-w)\Sigma=w+(1-2w)\Sigma$.

If $w<\frac{1}{2}$, the minimum of $L_w$ is equal to $w$, achieved when $\Sigma=0$,
that is, when $r\equiv 1$; if $w>\frac{1}{2}$, the minimum of $L_w$ is equal to 
$1-w$, achieved when $\Sigma=1$, that is, when $r\equiv 0$; if $w=\frac{1}{2}$,
the function $L_w$ is constant and equal to $\frac{1}{2}$, that is, all rules
are equally good and we might as well toss a coin.

$\bullet$ \emph{Proof that $(\rho,\alpha)\neq(\rho',\alpha')$ implies 
$G_{\rho,\alpha}\neq G_{\rho',\alpha'}$ (Section \ref{sec:Main}):}

Since $G'_{\rho,\alpha}(1)=-2\rho$, two different $\rho$ yield for sure two different
functions. Suppose $\rho=\rho'$ and $\alpha'>\alpha\ge 0$. 
Then, $\lim_{\eta\to\infty}\eta^{\rho-\alpha}G_{\rho,\alpha}(\eta)<\infty$
and $\lim_{\eta\to\infty}\eta^{\rho-\alpha}G_{\rho,\alpha'}(\eta)=\infty$,
hence $G_{\rho,\alpha}$ and $G_{\rho,\alpha'}$ must be different.

$\bullet$ \emph{Proof of Hallam's lemma (Lemma \ref{lem:preorder2poset}):}

This proof can be found in Hallam's PhD thesis \cite{Hallam2015}; we include it here
for the reader's convenience, and because his statement does not correspond 
completely to the proof. 

Suppose $\bar x\preceq\bar y$ and $\bar y\preceq\bar x$. Then $x_1\le y_1$ for
some $x_1\in\bar x$ and $y_1\in\bar y$. But $\bar y\preceq\bar x$ implies
there must be some $x_2\in\bar x$ such that $y_1\le x_2$. But then, there must be
some $y_2\in\bar y$ such that $x_2\le y_2$. And so on.
At some point we must have an equality of elements, since the set is finite.
Then the classes $\bar x$ and $\bar y$ must coincide.

$\bullet$ \emph{Proof that the hypothesis of Hallam's lemma holds in our case:}

We can assume everywhere that $y\ge z$. Otherwise, at any time we can interchange the roles
of $y$ and $z$. Assume that $(\rho_0,\alpha_0)\le(\rho_1,\alpha_1)$, with elements
$(x_0,y_0,z_0,t_0)\in(\rho_0,\alpha_0)$ and 
$(x_1,y_1,z_1,t_1)\in(\rho_1,\alpha_1)$, such that
$(x_0,y_0,z_0,t_0)\le(x_1,y_1,z_1,t_1)$, and assume in addition that
they belong to the transitive reduction in $(\mathbb T,\le)$. 
It is easy to find that there are two cases: 
Either (I) $\rho_1=\rho_0+1$ and $\alpha_1=\alpha_0-1$, or
(II) $\rho_1=\rho_0+1$ and $\alpha_1=\alpha_0+1$.

Take a generic element 
$(x,y,z,t)\in(\rho_0,\alpha_0)$.
In case (I), we have in particular $\alpha_0>0$, which implies $y>0$. 
Then  
$(x,y,z,t)\le (x+1,y-1,z,t)\in(\rho_1,\alpha_1)$.
In case (II), if $t>0$, take 
$(x,y+1,z,t-1)\in(\rho_1,\alpha_1)$,
and if $t=0$, take
$(x+1,y,z-1,t)\in(\rho_1,\alpha_1)$.
Note that $z=t=0$ cannot happen.

If the elements of $(\rho_0,\alpha_0)$
and $(\rho_1,\alpha_1)$ are not related by the transitive reduction, the 
argument can be iterated through a chain of elements related by the 
transitive reduction. Note that we have also deduced en passant the
transitive reduction of the quotient poset.

$\bullet$ \emph{Computation of the number of voting tables (Section \ref{sec:ComputSoft}):}

The \emph{Whitney numbers} $W_\rho$ of
a finite ranked poset are defined as the number of elements in rank level $\rho$. 
We compute the total number of tables in 
${\mathbb T}$ by computing first the Whitney numbers. We assume transposed 
tables are identified.

For odd positive ranks $\rho$, the possible
values of the pair $(x,t)$ are $(\rho+r,r)$, for $r=0,1,\dots,\frac{n-\rho}{2}$.

For each fixed $r$, the possible values of the pair $(y,z)$, with $y\ge z$,
are $(n-\rho-2r-s,s)$, for $s=0,1,\frac{n-\rho-2r}{2}$.
Thus, there are $\frac{n-\rho-2r}{2}+1$ such pairs. 
adding up these quantities from $r=0$ to $\frac{n-\rho}{2}$, yields
$\frac{1}{8}(n-\rho+4)(n-\rho+2)$.

A similar counting gives $\frac{1}{8}(n-\rho+3)(n-\rho+1)$ for even non-negative ranks $\rho$, which is the
same number as the odd rank immediately above. The case of negative ranks is deduced by symmetry.
Therefore we can write
\begin{equation*}
W_\rho = 
\begin{cases}
\frac{1}{8}(n-\rho+4)(n-\rho+2), & \text{if $\rho=1,3,\dots,n$} \\
W_{\rho+1}, & \text{if $\rho=0,2,\dots,n-1$} \\
W_{-\rho}, & \text{if $-n\le \rho\le -1$} 
\end{cases}
\end{equation*}
A simple but tedious computation, adding up all the Whitney numbers for $-n\le \rho \le n$, yields
the total number of tables $\sum_{\rho=-n}^n W_\rho = \frac{1}{24}(2n^3+15n^2+34n+21)$ \ .

The number $\max_{\rho} W_{\rho}$ is the maximal cardinality of an antichain: Indeed, since 
there is a unique minimal and a unique maximal element in $({\mathbb T},\le)$, any 
element of the poset is comparable to the minimal
and maximal element of the poset, and therefore to some element of
the most populated rank level. 
This implies that there cannot be more that $\max_{\rho} W_\rho$ elements
in any antichain. Ranks $-1$, $0$, and $1$ are the most populated, and its
Whitney number is $\frac{1}{8}(n+3)(n+1)$.

For completeness, let us just mention that the number of tables 
in the original poset $(\mathbb T,\le)$, without identifying transposed tables, is 
$\frac{1}{24}(4n^3+24n^2+44n+24)$, and the most populated rank has 
$\frac{1}{4}(n+3)(n+1)$ elements. 
In the quotient poset of $(\rho,\alpha)$-tables, there are $\frac{1}{2}(n+2)(n+1)$ tables
and a maximum of $\frac{1}{2}(n+1)$ members in any antichain. 
All computations are similar to those shown here.

\bibliographystyle{abbrv}
\bibliography{GenDecisionRules}

\begin{thebibliography}{10}

\bibitem{AlabertFarre2022}
A.~Alabert and M.~Farr{\'e}.
\newblock The doctrinal paradox: comparison of decision rules in a
  probabilistic framework.
\newblock {\em Social Choice and Welfare}, 58(4):863--895, 2022.

\bibitem{Berger1985}
J.~O. Berger.
\newblock {\em {Statistical decision theory and Bayesian analysis; 2nd ed.}}
\newblock Springer {S}eries in {S}tatistics. Springer, New York, 1985.

\bibitem{10.2307/2348873}
P.~J. Boland.
\newblock Majority systems and the {C}ondorcet {J}ury {T}heorem.
\newblock {\em Journal of the Royal Statistical Society. Series D (The
  Statistician)}, 38(3):181--189, 1989.

\bibitem{Bovens2006}
L.~Bovens and W.~Rabinowicz.
\newblock Democratic answers to complex questions -- an epistemic perspective.
\newblock {\em Synthese}, 150(1):131--153, 2006.

\bibitem{Bozbay2019}
{\.I}.~Bozbay.
\newblock Truth-tracking judgment aggregation over interconnected issues.
\newblock {\em Social Choice and Welfare}, 53(2):337--370, 2019.

\bibitem{BozbayDietrichPeters2014}
{\.I}.~Bozbay, F.~Dietrich, and H.~Peters.
\newblock Judgment aggregation in search for the truth.
\newblock {\em Games and Economic Behavior}, 87:571--590, 2014.

\bibitem{deClippel201534}
G.~de~Clippel and K.~Eliaz.
\newblock Premise-based versus outcome-based information aggregation.
\newblock {\em Games and Economic Behavior}, 89:34--42, 2015.

\bibitem{DIETRICH2006286}
F.~Dietrich.
\newblock Judgment aggregation: (im)possibility theorems.
\newblock {\em Journal of Economic Theory}, 126(1):286--298, 2006.

\bibitem{DietrichList2007}
F.~Dietrich and C.~List.
\newblock Judgment aggregation by quota rules: Majority voting generalized.
\newblock {\em Journal of Theoretical Politics}, 19, 2007.

\bibitem{DietrichList2008}
F.~Dietrich and C.~List.
\newblock Judgment aggregation without full rationality.
\newblock {\em Social Choice and Welfare}, 31(1):15--39, 2008.

\bibitem{dietrich_spiekermann_2013}
F.~Dietrich and K.~Spiekermann.
\newblock Epistemic democracy with defensible premises.
\newblock {\em Economics and Philosophy}, 29(1):87–120, 2013.

\bibitem{10.1093/mind/fzt074}
F.~Dietrich and K.~Spiekermann.
\newblock {Independent Opinions? On the Causal Foundations of Belief Formation
  and Jury Theorems}.
\newblock {\em Mind}, 122(487):655--685, 2013.

\bibitem{DokowHolzman2009}
E.~Dokow and R.~Holzman.
\newblock Aggregation of binary evaluations for truth-functional agendas.
\newblock {\em Social Choice and Welfare}, 32(2):221--241, 2009.

\bibitem{Fallis2005}
D.~Fallis.
\newblock Epistemic value theory and judgment aggregation.
\newblock {\em Episteme}, 2(1):39–55, 2005.

\bibitem{Fawcett:2006:IRA:1159473.1159475}
T.~Fawcett.
\newblock An introduction to {ROC} analysis.
\newblock {\em Pattern Recognition Letters}, 27(8):861--874, 2006.

\bibitem{Gardenfors2006}
P.~G{ä}rdenfors.
\newblock A representation theorem for voting with logical consequences.
\newblock {\em Economics and Philosophy}, 22(2):181–190, 2006.

\bibitem{Grofman1983}
B.~Grofman, G.~Owen, and S.~L. Feld.
\newblock Thirteen theorems in search of the truth.
\newblock {\em Theory and Decision}, 15(3):261--278, 1983.

\bibitem{Hallam2015}
J.~Hallam.
\newblock {\em Quotient Posets and the Characteristic Polynomial}.
\newblock PhD thesis, Michigan State University, 2015.

\bibitem{Hartmann-al2010}
S.~Hartmann, G.~Pigozzi, and J.~Sprenger.
\newblock {Reliable Methods of Judgement Aggregation}.
\newblock {\em Journal of Logic and Computation}, 20(2):603--617, 2010.

\bibitem{Kornhauser1992169}
L.~A. Kornhauser.
\newblock Modeling collegial courts {I}: Path-dependence.
\newblock {\em International Review of Law and Economics}, 12(2):169 -- 185,
  1992.

\bibitem{Kornhauser1992b}
L.~A. Kornhauser.
\newblock Modeling collegial courts, {II}. {L}egal doctrine.
\newblock {\em Journal of Law, Economics, and Organization}, 8(3):441--470,
  1992.

\bibitem{KornhauserSager1993}
L.~A. Kornhauser and L.~G. Sager.
\newblock The one and the many: Adjudication in collegial courts.
\newblock {\em California Law Review}, 81(1):1--51, 1993.

\bibitem{Lang2017}
J.~Lang, G.~Pigozzi, M.~Slavkovik, L.~van~der Torre, and S.~Vesic.
\newblock A partial taxonomy of judgment aggregation rules and their
  properties.
\newblock {\em Social Choice and Welfare}, 48(2):327--356, 2017.

\bibitem{List2005}
C.~List.
\newblock The probability of inconsistencies in complex collective decisions.
\newblock {\em Social Choice and Welfare}, 24(1):3--32, 2005.

\bibitem{List2012}
C.~List.
\newblock The theory of judgment aggregation: an introductory review.
\newblock {\em Synthese}, 187(1):179--207, 2012.

\bibitem{list_pettit_2002}
C.~List and P.~Pettit.
\newblock Aggregating sets of judgments: An impossibility result.
\newblock {\em Economics and Philosophy}, 18(1):89–110, 2002.

\bibitem{ListPettit2004}
C.~List and P.~Pettit.
\newblock Aggregating sets of judgments: Two impossibility results compared.
\newblock {\em Synthese}, 140(1):207--235, 2004.

\bibitem{ListPolak2010441}
C.~List and B.~Polak.
\newblock Introduction to judgment aggregation.
\newblock {\em Journal of Economic Theory}, 145(2):441--466, 2010.
\newblock Judgment Aggregation.

\bibitem{List2009-LISJAA}
C.~List and C.~Puppe.
\newblock Judgment aggregation.
\newblock In P.~Anand, P.~Pattanaik, and C.~Puppe, editors, {\em Handbook of
  Rational and Social Choice}. Oxford University Press, 2009.

\bibitem{MillerOsherson2009}
M.~K. Miller and D.~Osherson.
\newblock Methods for distance-based judgment aggregation.
\newblock {\em Social Choice and Welfare}, 32(4):575--601, 2009.

\bibitem{Miyashita2021}
M.~Miyashita.
\newblock Premise-based vs conclusion-based collective choice.
\newblock {\em Social Choice and Welfare}, 57(2):361--385, 2021.

\bibitem{Mongin2012-MONTDP}
P.~Mongin.
\newblock The doctrinal paradox, the discursive dilemma, and logical
  aggregation theory.
\newblock {\em Theory and Decision}, 73(3):315--355, 2012.

\bibitem{NEHRING20111488}
K.~Nehring and M.~Pivato.
\newblock Incoherent majorities: The {M}c{G}arvey problem in judgement
  aggregation.
\newblock {\em Discrete Applied Mathematics}, 159(15):1488--1507, 2011.

\bibitem{NehringPuppe2008}
K.~Nehring and C.~Puppe.
\newblock Consistent judgement aggregation: the truth-functional case.
\newblock {\em Social Choice and Welfare}, 31(1):41--57, 2008.

\bibitem{oeis-org}
{OEIS Foundation}.
\newblock The on-line encyclopedia of integer sequences, sequence number
  {A}000372.
\newblock \url{https://oeis.org/}.

\bibitem{PaulyVanHess2006}
M.~Pauly and M.~van Hees.
\newblock Logical constraints on judgement aggregation.
\newblock {\em Journal of Philosophical Logic}, 35(6):569--585, 2006.

\bibitem{Peleg2012}
B.~Peleg and S.~Zamir.
\newblock Extending the {C}ondorcet jury theorem to a general dependent jury.
\newblock {\em Social Choice and Welfare}, 39(1):91--125, 2012.

\bibitem{Pettit2001}
P.~Pettit.
\newblock Deliberative democracy and the discursive dilemma.
\newblock {\em Philosophical Issues}, 11(1):268--299, 2001.

\bibitem{Pigozzi2006}
G.~Pigozzi.
\newblock Belief merging and the discursive dilemma: an argument-based account
  to paradoxes of judgment aggregation.
\newblock {\em Synthese}, 152(2):285--298, 2006.

\bibitem{Pigozzi2009}
G.~Pigozzi, M.~Slavkovik, and L.~van~der Torre.
\newblock A complete conclusion-based procedure for judgment aggregation.
\newblock In F.~Rossi and A.~Tsoukias, editors, {\em Algorithmic Decision
  Theory}, pages 1--13. Springer Berlin Heidelberg, 2009.

\bibitem{Pivato2013}
M.~Pivato.
\newblock Voting rules as statistical estimators.
\newblock {\em Social Choice and Welfare}, 40(2):581--630, 2013.

\bibitem{PIVATO201751}
M.~Pivato.
\newblock Epistemic democracy with correlated voters.
\newblock {\em Journal of Mathematical Economics}, 72:51--69, 2017.

\bibitem{TerzopoulouEndrissSAGT2019}
Z.~Terzopoulou and U.~Endriss.
\newblock Optimal truth-tracking rules for the aggregation of incomplete
  judgments.
\newblock In {\em Proceedings of the 12th International Symposium on
  Algorithmic Game Theory (SAGT-2019)}, 2019.

\end{thebibliography}
\end{document}